\documentclass{article}
\usepackage{amsmath}
\usepackage{amsfonts}
\usepackage{amssymb}
\usepackage{graphicx}
\usepackage{algpseudocode}
\usepackage{algorithm}
\usepackage[utf8]{inputenc}
\usepackage[dvipsnames]{xcolor}

\renewcommand{\Re}{\mathbb{R}}

\newcommand{\jhat}{\hat{\textit{\j}}}
\newcommand{\hhat}{\hat{\textit{h}}}

\newcommand{\ihat}{\hat{\textit{\i}}}

\newtheorem{proposition}{Proposition}[section]

\newtheorem{remark}{Remark}[section]

\title{Avoiding local minima in multilayer network optimization by incremental training}
\author{Alberto De Santis\thanks{Dipartimento di Ingegneria informatica automatica e gestionale ``A.Ruberti" -- ``Sapienza" Università di Roma. e-mails: \texttt{desantis@diag.uniroma1.it}, \texttt{liuzzi@diag.uniroma1.it},
\texttt{lucidi@diag.uniroma1.it},
\texttt{tronci@diag.uniroma1.it}
}, Giampaolo Liuzzi$^*$,\\ Stefano Lucidi$^*$, Edoardo M. Tronci$^*$}
\date{March 2021}

\begin{document}

\maketitle

\begin{quote}
    \small
    {\bf Abstract}. Training a large multilayer neural network can present many difficulties due to the large number of useless stationary points. These points usually attract the minimization algorithm used during the training phase, which therefore results inefficient. Extending some results proposed in literature for shallow networks, we propose the mathematical characterization of a class of such stationary points that arise in deep neural networks training. Availing such a description, we are able to define an incremental training algorithm that avoids getting stuck in the region of attraction of these undesirable stationary points. 
\end{quote}

\begin{quote}
\end{quote}

\section{Introduction}
Training any kind of neural networks is doubtless an extremely difficult task because in general it requires to minimize a non-convex objective loss function
\begin{equation*}
    \min_{x \in \mathbb{R}^n} \ f(x),
\end{equation*}
which may depend on a considerable number of training parameters. 
Another crucial issue regarding this problem is doubtless the presence of a consistent number of ``useless" stationary points, i.e. stationary points that are quite far from those relevant minimum points of the loss function. Such stationary points could be in principle the cause of the so called plateaus phenomenon which plagues the minimization of $f(x)$. 

\par\smallskip

Many interesting mathematical  characterizations of such stationary points have been proposed in literature (\cite{baldi:1989,baldi:2012,brady:1989,gori:1992,sontag:1989}). For instance, in \cite{fukumizu:2000} it has been theoretically proved  that even shallow networks (i.e. networks with only one single hidden layer) usually have many of such undesirable stationary points.

\par\medskip

An interesting topic of research consists in studying how to fruitfully use the mathematical characterization of stationary points to develop ad-hoc optimization algorithms for neural networks training.
In  \cite{fukumizu:2000}, it has been proved  that a subset of critical points of a shallow network with $H-1$ neurons in the hidden layer gives rise to submanifolds of critical points for a larger network with $H$ neurons.  

\par\smallskip

Our aim in the paper is twofold. First  a more theoretical contribution is given by extending the result of \cite{fukumizu:2000} to deep multilayer neural networks, i.e. networks with more than one layer. More precisely, we show that classes of stationary points of a given network derive from stationary points of smaller size networks, which are obtained from the larger size one by discarding an arbitrary (possibly large) number of neurons. On the other hand, the structure of such manifolds of stationary points shows that their number grows exponentially with the dimension of the network.

We then provide a computational contribution which exploits the extended characterization of classes of uninteresting stationary points to define a new optimization strategy which avoids getting stuck in their regions of attraction. More specifically, the proposed strategy is based on an incremental approach that produces points which have a loss function value smaller than that associated to the  particular classes of those useless stationary points. 

The computational experience shows that the proposed algorithm is able to effectively take advantage of the above mentioned mathematical characterization of stationary points on a variety of learning problems.

% The numerical experience carried the proposed incremental training algorithm is able to avoid plateau regions and to attain better loss function values. We propose a simple but effective incremental training algorithm (ITA), which differs from the approaches previously mentioned by being an alternative way of training a shallow neural network without modifying its structure or altering the loss landscape in any way. 
% The authors in \cite{fukumizu:2000} have proved for shallow networks that a subset of critical points corresponding to the global minimum of a smaller network with $H-1$ neurons gives rise to many local minima and saddle points of a larger network with $H$ neurons. The idea behind our incremental approach is to extend the work in \cite{fukumizu:2000} by proving for shallow networks that every stationary point for a network with $H$ neurons corresponds to stationary points for a network with $H+K$ neurons. 
% This result allows to train a shallow neural network of a given dimension starting from a smaller network and recursively increasing the number of neurons. At every iteration, the training of the larger network is performed by properly choosing both the starting point and the minimization technique. 
% By carefully choosing the initial parameters configuration, this method allows to progressively enlarge the neural network dimension during the training phase, which definitely aids the optimizer to avoid plateaus and speed up the training phase due to the fact of acting in a lower loss surface dimension with respect to the original problem.

The remainder of the paper is organized as follows. We first introduce in Section \ref{sec:2} the notation used in this paper to define a general structure of our neural network. In Section \ref{sec:3} we consider the case of networks with more than one hidden layer and more outputs and when an arbitrary number of neurons are added to a given layer in the network. In Section \ref{sec:4} we further generalize the result  by allowing to add any number of neurons on any intermediate layer.  In Section \ref{sec:5} we formally state our incremental training algorithm (ITA) and in Section \ref{sec:6} we report a detailed numerical experience of our method over a significant set of test problems. We finally,  in Section \ref{sec:7} draw some conclusions and outline possible future developments.

\section{Topology and Notation}\label{sec:2}\noindent
Given a supervised learning problem of the form $\{(x^p, y^p)\}_{p=1}^{P}$, where $x^p \in \mathbb{R}^n$, $y^p \in \mathbb{R}^m$ and $P$ is the number of training samples, we know from Hornik (\cite{hornik:1991}, \cite{hornik:1989}, \cite{hornik:1990}) that any continuous function on a closed and bounded subset of $\mathbb{R}^n$ can be universally approximated by a multilayer neural network, by finding a particular set of parameters $\theta^*$ that result in the best function approximation $f(x, \theta^*)$ within a given function space $\mathcal{F} = \{f(x, \theta): \mathbb{R}^n \rightarrow \mathbb{R}^m| \theta \in \Theta\}$, where $\Theta$ represents the parameter space. A neural network can be seen as an acyclic oriented graph, as illustrated in Fig. \ref{fig1}:
\begin{figure}[ht]
    \centering
    \caption{General structure of a multilayer neural network}
    \includegraphics[width=\textwidth]{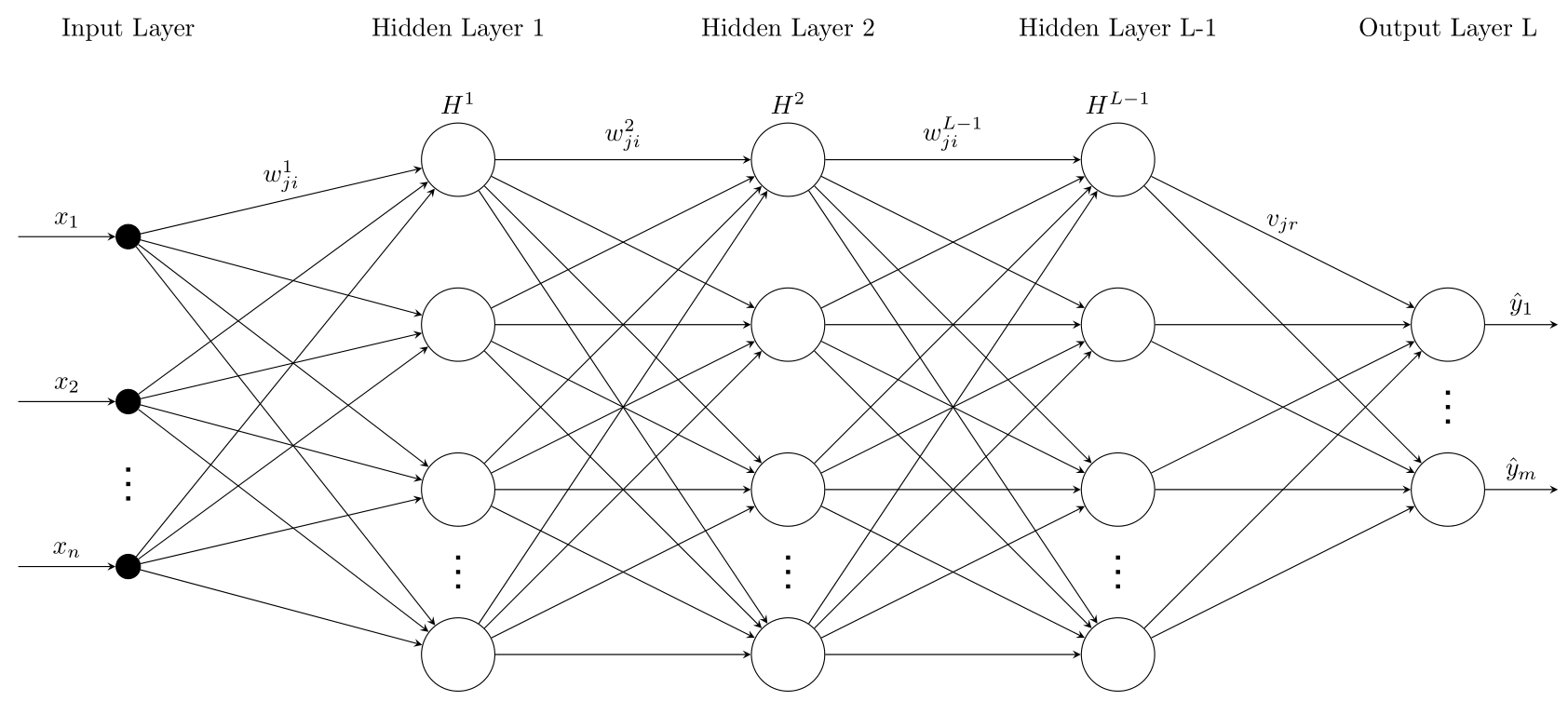}
    \label{fig1}
\end{figure}
the first layer holds the set of $n$ input nodes that connect each input component $x_{i} \in \mathbb{R},$ $i= 1, \ldots ,n$ to the network. Thus, $x^p = (x_{1}^{p}, \ldots, x_{n}^{p}), \hspace{0.2cm} p= 1, \ldots, P$ represents the input training set. A set of artificial neurons are distributed in $L$ layers, which are connected together in a chain such that the first $L-1$ hidden layers have no direct connections with the output. 
The output layer $L$ holds $m$ neural units, one for each output dimension $y_{r} \in \mathbb{R},$ $r= 1, \ldots, m$. Thus, $y^p = (y_{1}^{p}, \ldots, y_{m}^{p}),  \hspace{0.2cm} p= 1, \ldots, P$ represents the labels training set.
The overall length of the chain gives the depth of the model. Instead, the dimensionality of the hidden layers determines the width of the model.
We will only consider fully connected (dense) networks, i.e. each neuron of each layer is connected with all the neurons of the previous and the subsequent layers. We have assumed that exist neither connections between neurons of the same layer, nor feedback connections between outputs of a layer and inputs of the preceding layers. The activation function is the same for each neuron of every hidden layer. Furthermore, we have considered a bias term associated to each neuron unit included the output units.
To simplify the analysis, we have denoted with $\ell = 1, \ldots, L$ the index of layers and with $H_{\ell}$ the number of neurons in the $\ell$-th layer, where $H_{0} = n$ and $H_{L} = m$.
Let $\theta \in \mathbb{R}^q$ be the the vector that encases all the network parameters, where $q$ is the total number of network parameters, i.e.
\[
  q = \sum_{i=0}^{L-1} H_{i}\times H_{i+1} + \sum_{i=1}^L H_i.
\]
Thus, for every $\ell = 1, \ldots, L$ we can decompose it as 
\begin{eqnarray*}
&& (\theta^1)^T=(\sigma_1^1,(w_1^1)^T,\ldots,\sigma_{H_1}^1,(w_{H_1}^1)^T)\\
&& (\theta^2)^T=((\theta^1)^T,\sigma_1^2,(w_1^2)^T,\ldots,\sigma_{H_2}^2,(w_{H_2}^2)^T)\\
&& \ldots\qquad\ldots\qquad\ldots\qquad\ldots\qquad\ldots\\
&& \ldots\qquad\ldots\qquad\ldots\qquad\ldots\qquad\ldots\\
&& (\theta^L)^T=((\theta^{L-1})^T,\sigma_1^L,(w_1^L)^T,\ldots,\sigma_{H_L}^L,(w_{H_L}^L)^T)
\end{eqnarray*}
where 
$$\sigma_{j}^{\ell} \in \mathbb{R},\hspace{0.5cm} w_{j}^{\ell}\in \mathbb{R}^{H_{\ell-1}},\hspace{0.3cm} j=1,\ldots, H_{\ell}$$
are the bias term of neuron $j$ in layer $\ell$ and the weights vector of layer $\ell$ that assigns a scalar $w_{ji}^{\ell} \in \mathbb{R}$ to each arc, which is the weight between neuron $j$ of layer $\ell$ and neuron $i$ of layer $\ell-1$. %\red{ dovremmo indicare chi è $\theta$ e dire anche che $\Theta = R^{n. totale di parametri}$}
Note that $\theta^L = \theta$.\\
Being $g$ : $ R \rightarrow R$ the activation function, and denoting with $a_{j}^{\ell}$ the input of neuron $j$ in layer $\ell$, we obtain for neuron $j$ of the first layer
\begin{equation*}
    a_{j}^{1}(x, \theta^1) = \sum_{i=1}^{n} w_{ji}^{1}\; x_i + \sigma_{j}^{1}, \hspace{0.3cm}  j=1, \ldots, H_{1}
\end{equation*}
and for neuron $j$ of layer $\ell > 1$
\begin{equation*}
    a_{j}^{\ell}(x, \theta^{\ell}) = \sum_{i=1}^{H_{\ell-1}} w_{ji}^{\ell}\; g(a_{i}^{\ell - 1}(x, \theta^{\ell - 1})) + \sigma_{j}^{\ell}, \hspace{0.3cm}  j=1, \ldots, H_{\ell}
\end{equation*}
We further assume that the output units are linear. Thus, for $\ell=L$ we have $f_r(x, \theta^{L}) = a_{r}^{L}(x, \theta^{L}), \hspace{0.2cm} r= 1, \ldots, m$.\\
Under the assumptions stated, we can define an input-output mapping of the form
\begin{equation}
f_r(x, \theta^{L}) = \sum_{i=1}^{H_{L-1}} w_{ri}^{L}\; g(a_{i}^{L - 1}(x, \theta^{L - 1})) + \sigma_{r}^{L}, \hspace{0.3cm}  r=1, \ldots, m.
\label{eq1} 
\end{equation}
In order to choose an approximate function among all the possible functions in $\mathcal{F}$, we introduce the empirical risk 
\begin{equation}
R_{emp}(\theta)=\frac{1}{P}\sum_{p=1}^{P} \mathcal{L}_{p},
\label{eq2} 
\end{equation}
where $\mathcal{L}_{p} = \mathcal{L} (y^{p}, f(x^{p}, \theta^{L}))\geq 0$  evaluates the distance between the experimental data $y^p$ and the output generated by the model $f(x^{p}, \theta^{L})$. As known in the literature, $\mathcal{L}$ can assume different forms depending on the problem faced. The results in this paper are independent of the choice of the loss function.

The next step is to compute the derivative of $R_{emp}$ with respect to the network parameters $\theta$. In this regard, various methods have been proposed in literature, e.g. the back propagation approach. Here, we adopt the forward one since it is more suited for our needs. 

For every ${\ell} =1,\dots,L$, let $j=1,\ldots,H_\ell$  and $i=1,\ldots,H_{\ell-1}$. Denote by $\delta_{cj}$ the kronecker delta symbol. For the input layer ${\ell} =1$ we can write
\begin{eqnarray}
\nonumber	
	%&&\hbox{for}\quad c=1,\ldots,H_\ell\\ [0.3truecm]
%\nonumber	
	%&&\hbox{if }\quad {\ell} =1	\\ [0.3truecm]
		\label{eq3}
	&&\quad \frac{\partial a^{1}_c(x^{p}, \theta^1)}{\partial \sigma^{{1}}_{j}}=\delta_{cj},\\
	[0.3truecm]
	\label{eq4}
	&&\quad \frac{\partial a^{1}_c(x^{p}, \theta^1)}{\partial w^{{1}}_{j,i}}=	x_i \delta_{cj},\quad c=1,\ldots,H_\ell.
	\\ [0.3truecm]\nonumber
	\end{eqnarray}
	For any successive layer ${\ell} =2,\ldots,L$, we have that
	\begin{eqnarray}
	%&&\hbox{otherwise }	\\ [0.3truecm]	
	\label{eq5}
	&&\quad \frac{\partial a^{\ell}_c(x^{p}, \theta^\ell)}{\partial \sigma^{\ell}_{j}}=\delta_{c j},
	\\ [0.3 truecm]	
	\label{eq6}
	&&\quad\frac{\partial a^\ell_c(x^p,\theta^\ell)}{\partial w^\ell_{i j}} = 
	g\left( a^{\ell -1}_i(x^p,\theta^{\ell-1})\right) \delta_{c j}
	%&&\quad \frac{\partial a^\ell_c (x^p, \theta^\ell)}{\partial w^\ell_{j i}}=	 g\biggl(a_i^{\ell-1}(x^p,\theta^{\ell-1})\biggr)\delta_{c j},\quad c=1,\ldots,H_\ell$.
	\end{eqnarray}
	%\\ [0.3truecm]
	%\nonumber
	%&&\hbox{where}\ \delta_{cj}\ \hbox{is the kronecker delta.} \\ [0.3truecm]	
	The previous relations evaluate the derivative of the node inputs with respects to their parameters. We now compute how the input parameters of the nodes of any layer influence the outputs of the nodes of any successive layer. For any $\ell<L$, let $q =\ell+1,\dots,L$ and $c=1,\ldots,H_q$. We have 
	\begin{eqnarray}
	%&&\hbox{If} \quad \ell<L\quad   \hbox{then for } \quad {q} =\ell+1,\dots,L, \quad c=1,\ldots,H_q	\\ [0.3truecm]
	\label{eq7}
	&&\quad\quad \frac{\partial a^{q}_c(x^{p}, \theta^{q}))}{\partial \sigma^{{\ell}}_{j}}=	\sum_{h=1}^{H_{q-1}}w^q_{ch} g^\prime\left(a_{h}^{q-1}(x^{p},\theta^{q-1})\right) \frac{\partial a_{h}^{q-1}(x^{p},\theta^{q-1})}{\partial \sigma^{{\ell}}_{j}},\\ [0.3truecm]
	\label{eq8}
	&&\quad\quad \frac{\partial a^{q}_c(x^{p}, \theta^q))}{\partial w^{{\ell}}_{ji}}=	\sum_{h=1}^{H_{q-1}}w^q_{ch} g^\prime\left(a_{h}^{q-1}(x^{p},\theta^{q-1})\right) \frac{\partial a_{h}^{q-1}(x^{p},\theta^{q-1})}{\partial w^{{\ell}}_{ji}}.
\end{eqnarray}

Thus, for every\quad  ${\ell} =1,\dots,L$,\quad $j=1,\ldots,H_\ell$, \quad $i=1,\ldots,H_{\ell-1}$, we can write
\begin{eqnarray}
	\label{eq9}
	&&\frac{\partial R_{emp}(\theta^L)}{\partial \sigma^{{\ell}}_{j}}  = \frac{1}{P} \sum_{p=1}^{P} \sum_{r=1}^{m} { \mathcal{L}^\prime(y^{p}, f_r(x^{p}, \theta^L))}
	\frac{\partial a^{L}_r(x^{p}, \theta^L))}{\partial \sigma^{{\ell}}_{j}},\\ [0.3truecm]
	\label{eq10}
	&&\frac{\partial R_{emp}(\theta^L)}{\partial w^{{\ell}}_{ji}}  = \frac{1}{P} \sum_{p=1}^{P} \sum_{r=1}^{m}{ \mathcal{L}^\prime(y^{p}, f_r(x^{p}, \theta^L))} 
	\frac{\partial a^{L}_r(x^{p}, \theta^L))}{\partial w^{{\ell}}_{ji}}.
\end{eqnarray}
%------------------------------------

%\input{old_material.tex}

\section{Adding $K$ neurons to one of the layers}\label{sec:3}
In this section, we consider the case in which we add $K$ neurons to the $\hat\ell$-th layer of a given neural network. Here, we denote by $\theta$ the vector of parameters of the smaller network and by $\hat\theta$ the vector of the bigger network.

In order to compute the empirical risk of the new network, we need to define the following quantities. Starting from  
\begin{equation}
	\label{eq11}
	\hat a_j^{1 }(x,\hat \theta^1)=\sum_{i=1}^{n}
	\hat w_{ji}^{{1} }\;x_i+\hat \sigma^{1 }_{j},\qquad {j} =1,\dots,H_1,	
\end{equation}
for $\ell =2,\dots,\hat\ell-1$, we have	
\begin{equation}
	\label{eq12} 
	\qquad \hat a_j^{{\ell} }(x,\hat\theta^\ell)=\sum_{i=1}^{H_{\ell-1}} \hat w_{ji}^{{\ell} }\;{g(\hat a_i^{{\ell} -1}(x,\hat\theta^{{\ell} -1}))}+\hat \sigma^{{\ell} }_{j},\qquad {j} =1,\dots,H_\ell.
\end{equation}
Furthermore
\begin{eqnarray}
	\label{eq13}
	&&\hspace*{-0.7cm} \hat a_j^{{\hat\ell} }(x,\hat\theta^{\hat\ell})=\sum_{i=1}^{H_{{\hat\ell}-1}} \hat w_{ji}^{{\hat\ell} }\;{g(\hat a_i^{{{\hat\ell}} -1}(x,\hat\theta^{{\hat\ell} -1}))}+\hat \sigma^{{\hat\ell} }_{j},\qquad {j} =1,\dots,H_{\hat\ell},	\\ [0.3truecm]
	\label{eq14}
	&&\hspace*{-0.7cm} \hat a_{\jhat}^{{\hat\ell} }(x,\hat\theta^{\hat\ell})=\sum_{i=1}^{H_{{\hat\ell}-1}} \hat w_{{\jhat}i}^{{\hat\ell} }\;{g(\hat a_i^{{{\hat\ell}} -1}(x,\hat\theta^{{\hat\ell} -1}))}+\hat \sigma^{{\hat\ell} }_{{\jhat}},\quad {\jhat} =H_{\hat\ell}+1,\dots,{H_{\hat\ell}+K},	\\ [0.3truecm]
	\label{eq15}
	&&\hspace*{-0.7cm} \hat a_j^{{\hat\ell+1} }(x,\hat\theta^{\hat\ell+1})=\sum_{i=1}^{H_{{\hat\ell}}} \hat w_{ji}^{{\hat\ell}+1 }\;{g(\hat a_i^{{{\hat\ell}}}(x,\hat\theta^{{\hat\ell} }))}+\\
	&&\hspace*{-0.7cm} \qquad\quad \sum_{\ihat=H_{{\hat\ell}}+1}^{{H_{\hat\ell}+K}} \hat w_{j\ihat}^{{\hat\ell}+1 }\;{g(\hat a_{\ihat}^{{{\hat\ell}}}(x,\hat\theta^{{\hat\ell} }))}+\hat \sigma^{{\hat\ell+1} }_{j},\qquad {j} =1,\dots,H_{\hat\ell+1}.	\nonumber
\end{eqnarray}

For $\ell =\hat\ell+2,\dots,L$
\begin{equation}
	\label{eq16}
	\qquad \hat a_j^{{\ell} }(x,\hat\theta^\ell)=\sum_{i=1}^{H_{\ell-1}} \hat w_{ji}^{{\ell} }\;{g(\hat a_i^{{\ell} -1}(x,\hat\theta^{{\ell} -1}))}+\hat \sigma^{{\ell} }_{j},\qquad {j} =1,\dots,H_\ell,
\end{equation}
and
\begin{eqnarray}
	\label{eq17}
	&& \hat f_r^{}(x, \hat\theta^{L})=\hat a_r^{{L} }(x,\hat\theta^L), \qquad {r} =1,\dots,m.	
\end{eqnarray}
\par\bigskip\noindent

Therefore, the empirical risk of the new network can be written as
\begin{eqnarray}\label{eq18}
&&\hat R_{emp}(\hat\theta^{})=\frac{1}{P}\sum_{p=1}^{P} \mathcal{ \hat L} (y^{p}, \hat f^{}(x^{p}, {\hat\theta^{L}})),\end{eqnarray}
where $\hat f:\Re^n\to\Re^m$.
\par\medskip\noindent

The vector $\hat\theta$ can be obtained starting from $\theta$ according to some particular rules. In \cite{fukumizu:2000}, some of these rules have been defined in the case of neural networks with only one hidden layer, i.e. the mappings $\alpha$, $\beta$ and $\gamma$. Here, we extend those mappings to the more general case of deep neural networks, i.e. networks with more than one hidden layer, as shown in Fig. \ref{embeddings}:
\begin{figure}[ht]
    \centering
    \caption{Canonical embeddings of a multilayer neural network}
    \includegraphics[scale=0.18]{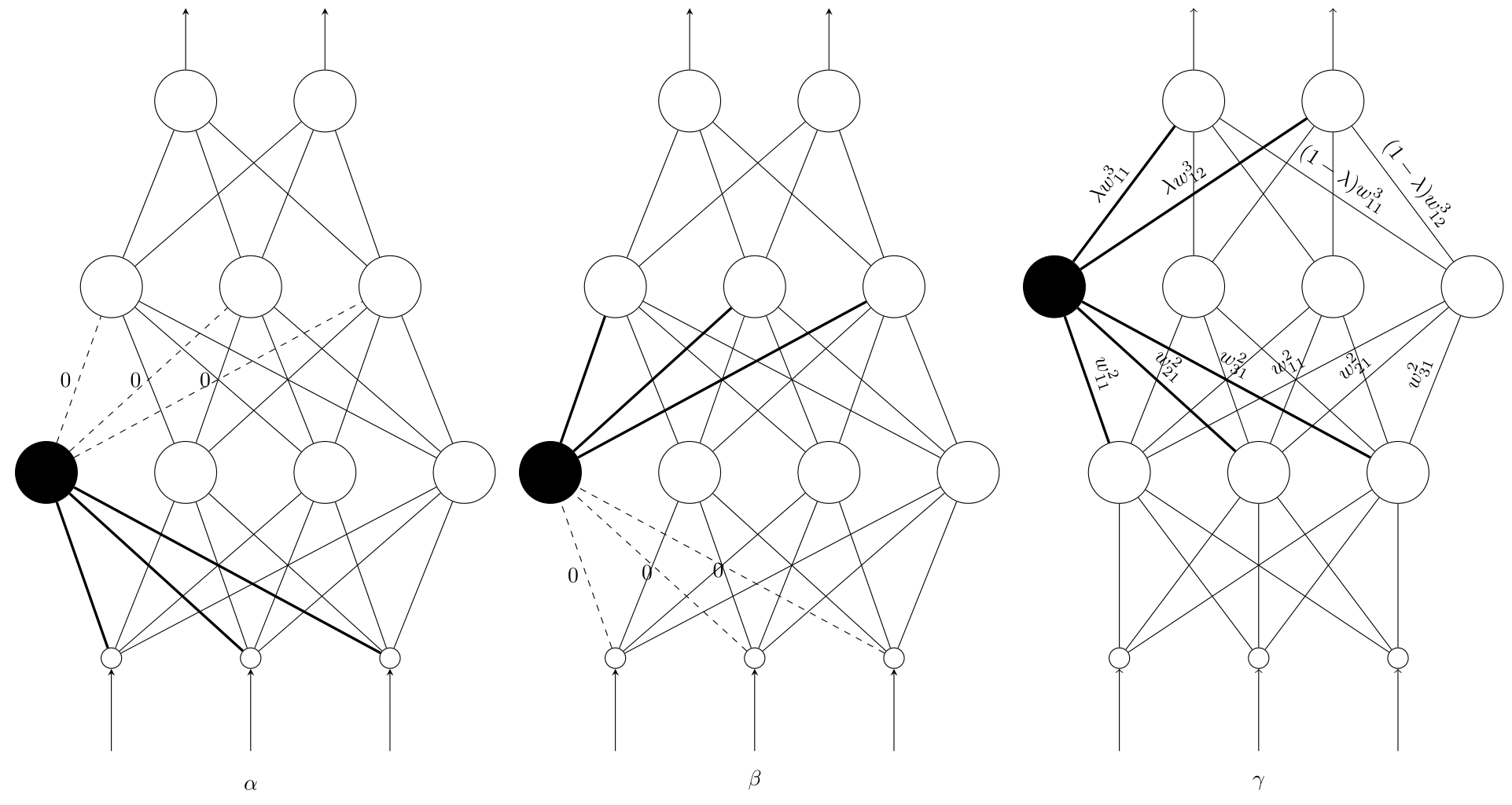}
    \label{embeddings}
\end{figure}
\begin{eqnarray*}
 \hat\theta=\alpha_{\zeta {v}}(\theta) = \biggl[& \hat \sigma_j^\ell=\sigma_j^\ell,\quad \hat w_{ji}^\ell=w_{ji}^\ell,\quad & \hbox{for } \left\{\begin{array}{l} {\ell} =1,\dots,L,\\ j=1,\ldots,H_\ell, \\ i=1,\ldots,H_{\ell-1}.\end{array}\right.\\
 &\hat \sigma_{\hat\jmath}^{\hat\ell}=\zeta_{\hat\jmath-H_{\hat\ell}},\quad  \hat w_{\jhat i}^{\hat\ell}=v_{\jhat-H_{\hat\ell} i},\quad & \hbox{for }\,\left\{\begin{array}{l}\jhat=H_{\hat\ell}+1,\ldots,H_{\hat\ell}+K, \\ i=1,\ldots,H_{{\hat\ell}-1}.\end{array}\right.\\
  &  \hat w_{j \ihat}^{\hat\ell+1}=0,\qquad & \hbox{for }\left\{\begin{array}{l} j=1,\ldots,H_{\hat\ell+1},\\ \ihat=H_{\hat\ell}+1,\ldots,H_{\hat\ell}+K\biggr],\end{array}\right.
\end{eqnarray*}
where $\zeta_{j}\in\Re$,  and $v_{j}\in\Re^{H_{\hat\ell - 1}}$, $j=1,\dots,K$.
\begin{eqnarray}
	\nonumber
	\hat\theta=\beta_{\zeta {s}}(\theta) = \Biggl[ \hat \sigma_j^\ell=\sigma_j^\ell,\quad \hat w_{ji}^\ell=w_{ji}^\ell,\quad && \hspace*{-0.2cm} \hbox{for }\left\{\begin{array}{l} {\ell} =1,\dots,L,\\ j=1,\ldots,H_\ell,\\ {\ell} \not=\hat\ell+1, \\ i=1,\ldots,H_{\ell-1}.\end{array}\right.\nonumber\\
	\label{eq19}
	\hat \sigma_{\jhat-H_{\hat\ell}}^{\hat\ell}=\zeta_{\jhat-H_{\hat\ell}},\quad\;  \hat w_{\jhat-H_{\hat\ell} i}^{\hat\ell}=0, && \hspace*{-0.2cm} \hbox{for }\left\{\begin{array}{l}  \jhat=H_{\hat\ell}+1,\ldots,H_{\hat\ell}+K, \\ i=1,\ldots,H_{{\hat\ell}-1}.\end{array}\right.\\
	\nonumber
	   \hat w_{j i}^{\hat\ell+1}=w_{j i}^{\hat\ell+1}, && \hspace*{-0.2cm} \hbox{for }\left\{\begin{array}{l}  j=1,\ldots,H_{\hat\ell+1},\\ i=1,\ldots,H_{\hat\ell}.\end{array}\right.\\
	\label{eq20}
	 \hat w_{j \ihat-H_{\hat\ell}}^{\hat\ell+1}=s_{j \ihat-H_{\hat\ell}}, && \hspace*{-0.2cm} \hbox{for }\left\{\begin{array}{l}  j=1,\ldots,H_{\hat\ell+1},\\ \ihat=H_{\hat\ell}+1,\ldots,H_{\hat\ell}+K.\end{array}\right.\\
	\nonumber
	  \hat \sigma_{j }^{\hat\ell+1}=\sigma_{j }^{\hat\ell+1}-\sum_{\ihat=H_{\hat\ell+1}}^{H_{\hat\ell}+K}s_{j \ihat-H_{\hat\ell}}\;g(\zeta_{\ihat - H_{\hat\ell}}), && \hspace*{-0.2cm} \hbox{for }\quad   j=1,\ldots,H_{\hat\ell+1}\Biggr],
\end{eqnarray}
where $\zeta_j\in\Re$, $s_j\in\Re^{H_{\hat\ell + 1}}$, $j=1,\dots,K$. Note that $\alpha_{\zeta\mathbf{0}}(\theta) = \beta_{\zeta\mathbf{0}}(\theta)$ trivially results. Thus, these two embeddings give the same critical point set. 
Letting $h\in \{1,\dots,H_{\hat\ell}\}$, we can then define the last embedding as

\begin{eqnarray}
	\nonumber
	\hat\theta=\gamma_{\lambda }(\theta) = \Biggl[ \hat \sigma_j^\ell=\sigma_j^\ell,\quad \hat w_{ji}^\ell=w_{ji}^\ell, && \hspace*{-0.2cm} \hbox{for }\left\{\begin{array}{l} {\ell} =1,\dots,L,\\ j=1,\ldots,H_\ell,\\ {\ell} \not=\hat\ell+1, \\ i=1,\ldots,H_{\ell-1}.\end{array}\right.\\
	\label{eq21}
	\qquad\qquad\qquad\;\hat \sigma_{\jhat}^{\hat\ell}=\sigma_{h}^{\hat\ell},\quad  \hat w_{\jhat i}^{\hat\ell}=w_{h i}^{\hat\ell}, && \hspace*{-0.2cm} \hbox{for }\left\{\begin{array}{l} \jhat=H_{\hat\ell}+1,\ldots,{H_{\hat\ell}+K}, \\ i=1,\ldots,H_{{\hat\ell}-1}.\end{array}\right.\\
	\label{eq22}
	\hat w_{j i}^{\hat\ell+1}=w_{j i}^{\hat\ell+1}, && \hspace*{-0.2cm} \hbox{for }\left\{\begin{array}{l} j=1,\ldots,H_{\hat\ell+1},\\ i=1,\ldots,H_{\hat\ell},\quad i\not=h.\end{array}\right.\\
	\label{eq23}
	 \hat w_{j h}^{\hat\ell+1}=\lambda_0 w_{j h}^{\hat\ell+1}, && \hspace*{-0.2cm} \hbox{for }\left\{\begin{array}{l}  j=1,\ldots,H_{\hat\ell+1}. \end{array}\right.\\
	\label{eq24}
	\hat w_{j \ihat}^{\hat\ell+1}=\lambda_{\ihat-H_{\hat\ell}} w_{j h}^{\hat\ell+1}, && \hspace*{-0.2cm} \hbox{for }\left\{\begin{array}{l}  j=1,\ldots,H_{\hat\ell+1},\\ \ihat=H_{\hat\ell}+1,\ldots,{H_{\hat\ell}+K}.\end{array}\right.\\
	\nonumber
	 \hat \sigma_{j }^{\hat\ell+1}=\sigma_{j }^{\hat\ell+1}, && \hspace*{-0.2cm} \hbox{for }\quad  j=1,\ldots,H_{\hat\ell+1}\Biggr],
\end{eqnarray}
for all  $\lambda_i \in \mathbb{R},\ i=0,\dots,K$ and such that $\sum_{i=0}^{K}\lambda_i = 1$.

\begin{proposition}\label{Prop1}
	For every  point  $\theta$, it results: 
\begin{eqnarray}\label{eq25}
&& \hat R_{emp}(\hat\theta)	=\hat R_{emp}(\alpha_{\zeta {v}}(\theta))= R_{emp}(\theta),\\
\label{eq26}
&& \hat R_{emp}(\hat\theta)	=\hat R_{emp}(\beta_{\zeta {s}}(\theta) )= R_{emp}(\theta),\\
\label{eq27}
&& \hat R_{emp}(\hat\theta)	=\hat R_{emp}(\gamma_{\lambda }(\theta))= R_{emp}(\theta).
\end{eqnarray}
\end{proposition}
\par\medskip\noindent
{\bf Proof}.
By the definitions of the three maps, we have for $\ell = 1$
\par\medskip\noindent
\begin{equation*}
	 \hat a_j^{1 }(x,\alpha_{\zeta {v}}(\theta)^1)=\hat a_j^{1 }(x,\beta_{\zeta {s}}(\theta)^1)=\hat a_j^{1 }(x,\gamma_{\lambda }(\theta)^1)= a_j^{1 }(x,\theta^1), \qquad {j} =1,\dots,H_1.
\end{equation*}	

Instead, for $\ell =2,\dots,\hat\ell-1$

\begin{equation*}
	\qquad \hat a_j^{{\ell} }(x,\alpha_{\zeta {v}}(\theta)^\ell)=\hat a_j^{{\ell} }(x,\beta_{\zeta {s}}(\theta)^\ell)=\hat a_j^{{\ell} }(x,\gamma_{\lambda }(\theta)^\ell)= a_j^{{\ell} }(x,\theta^\ell),\quad {j} =1,\dots,H_\ell,
\end{equation*}

\begin{equation}
	\label{eq28}
		\hat a_j^{{\hat\ell} }(x,\alpha_{\zeta {v}}(\theta)^{\hat\ell})=\hat a_j^{{\hat\ell} }(x,\beta_{\zeta {s}}(\theta)^{\hat\ell})=\hat a_j^{{\hat\ell} }(x,\gamma_{\lambda }(\theta)^{\hat\ell})= a_j^{{\hat\ell} }(x,\theta^{\hat\ell}),\quad {j} =1,\dots,H_{\hat\ell}.
\end{equation}
Now, let's consider the three diffent maps.
\par\medskip\noindent
If $\hat \theta= \alpha_{\zeta {v}}(\theta)$, we have that
\begin{eqnarray*}
	%\label{a-bar-41}
	&& \hat a_{\jhat}^{{\hat\ell} }(x,\hat\theta^{\hat\ell})=\sum_{i=1}^{H_{{\hat\ell}-1}}  v_{{\jhat}i}\;{g(\hat a_i^{{{\hat\ell}} -1}(x,\theta^{{\hat\ell} -1}))}+\zeta_{{\jhat}},\qquad {\jhat} =H_{\hat\ell}+1,\dots,{H_{\ell}+K},	\\ [0.3truecm]
	%\label{a-bar-5}
	&& \hat a_j^{{\hat\ell+1} }(x,\hat\theta^{\hat\ell+1})=\sum_{i=1}^{H_{{\hat\ell}}} \hat w_{ji}^{{\hat\ell}+1 }\;{g( \hat a_i^{{{\hat\ell}}}(x,\hat\theta^{{\hat\ell} }))}\\
	&&\hspace*{-0.3cm} \qquad\quad +\sum_{\ihat=H_{{\hat\ell}}+1}^{{H_{\ell}+K}}  0\;{g(\hat a_{\ihat}^{{{\hat\ell}}}(x,\hat\theta^{{\hat\ell} }))}+\hat \sigma^{{\hat\ell+1} }_{j}, \quad {j} =1,\dots,H_{\hat\ell+1}.	
\end{eqnarray*}
The definition of the map $\alpha_{\zeta {v}}$ and  (\ref{eq28}) imply that
\begin{eqnarray}\nonumber
	%\label{a-bar-51}
	&& \hat a_j^{{\hat\ell+1} }(x,\hat\theta^{\hat\ell+1})=\sum_{i=1}^{H_{{\hat\ell}}}  w_{ji}^{{\hat\ell}+1 }\;{g( a_i^{{{\hat\ell}}}(x,\theta^{{\hat\ell} }))}+\sigma^{{\hat\ell+1} }_{j}	\\ [0.3truecm]
	\nonumber
	&& \qquad\qquad\qquad\qquad= a_j^{{\hat\ell+1} }(x,\theta^{\hat\ell+1}),\quad {j} =1,\dots,H_{\hat\ell+1},		
\end{eqnarray}
which ensures that, with the definition of $\alpha_{\zeta {v}}$, and for $\ell =\hat\ell+2,\dots,L$
\begin{equation*}
	\qquad \hat a_j^{{\ell} }(x,\hat\theta^\ell)=\sum_{i=1}^{H_{\ell-1}}  w_{ji}^{{\ell} }\;{g( a_i^{{\ell} -1}(x,\theta^{{\ell} -1}))}+ \sigma^{{\ell} }_{j}= a_j^{{\ell} }(x,\theta^\ell),\qquad {j} =1,\dots,H_\ell, 
\end{equation*}
\begin{equation}
	\label{eq29}
	 \hat f_r(x,\hat\theta^{L})=f_r(x, \theta^{L}), \qquad {r} =1,\dots,m.
\end{equation}
In conclusion (\ref{eq2}), (\ref{eq18}) and (\ref{eq29}) prove (\ref{eq25}).
\par\medskip\noindent
If $\hat \theta= \beta_{\zeta {s}}(\theta)$, we have instead
\begin{eqnarray*}
	%\label{a-bar-41}
	&& \hat a_{\jhat}^{{\hat\ell} }(x,\hat\theta^{\hat\ell})=\sum_{i=1}^{H_{{\hat\ell}-1}}  0\;{g(\hat a_i^{{{\hat\ell}} -1}(x,\theta^{{\hat\ell} -1}))}+\zeta_{{\jhat}},\qquad {\jhat} =H_{\hat\ell}+1,\dots,{H_{\ell}+K},	\\ [0.3truecm]
	%\label{a-bar-5}
	&& \hat a_j^{{\hat\ell+1} }(x,\hat\theta^{\hat\ell+1})=\sum_{i=1}^{H_{{\hat\ell}}}  w_{ji}^{{\hat\ell}+1 }\;{g( \hat a_i^{{{\hat\ell}}}(x,\hat\theta^{{\hat\ell} }))}+\sum_{\ihat=H_{{\hat\ell}}+1}^{{H_{\ell}+K}}  s_{j\ihat}\;{g(\zeta_{\ihat})}+ \sigma^{{\hat\ell+1} }_{j}+\\ [0.3truecm]
	&&\qquad\qquad	-\sum_{\ihat=H_{{\hat\ell}}+1}^{{H_{\ell}+K}}  s_{j\ihat}\;{g(\zeta_{\ihat})}, \quad {j} =1,\dots,H_{\hat\ell+1}.
\end{eqnarray*}
Now, by using  (\ref{eq28}) we obtain
\begin{eqnarray}\nonumber
	%\label{a-bar-51}
	&& \hat a_j^{{\hat\ell+1} }(x,\hat\theta^{\hat\ell+1})=\sum_{i=1}^{H_{{\hat\ell}}}  w_{ji}^{{\hat\ell}+1 }\;{g( a_i^{{{\hat\ell}}}(x,\theta^{{\hat\ell} }))}+ \sigma^{{\hat\ell+1} }_{j}	\\ [0.3truecm]
	\nonumber
	&& \qquad\qquad\qquad\qquad= a_j^{{\hat\ell+1} }(x,\theta^{\hat\ell+1}),\quad {j} =1,\dots,H_{\hat\ell+1},
\end{eqnarray}

and for $\ell =\hat\ell+2,\dots,L$

\begin{eqnarray}
	\nonumber
	&&\qquad \hat a_j^{{\ell} }(x,\hat\theta^\ell)=\sum_{i=1}^{H_{\ell-1}}  w_{ji}^{{\ell} }\;{g( a_i^{{\ell} -1}(x,\theta^{{\ell} -1}))}+ \sigma^{{\ell} }_{j}	\\ [0.3truecm]
	\nonumber
	&&\qquad\qquad\qquad\qquad= a_j^{{\ell} }(x,\theta^\ell),\qquad {j} =1,\dots,H_\ell,
\end{eqnarray}
\[
	\hat f_r^{}(x, \hat\theta^{L})=f_r^{}(x, \theta^{L}), \qquad {r} =1,\dots,m,	
\]
thus proving (\ref{eq26}).
\par\medskip\noindent
Finally, if $\hat \theta= \gamma_{\lambda }(\theta) $ we have
\begin{eqnarray*}
	\label{a-bar-41}
	&& \hat a_{\jhat}^{{\hat\ell} }(x,\hat\theta^{\hat\ell})=\sum_{i=1}^{H_{{\hat\ell}-1}}  w_{h i}^{{\hat\ell} }\;{g(\hat a_i^{{{\hat\ell}} -1}(x,\theta^{{\hat\ell} -1}))}+\sigma_{{h}}^{\hat\ell}= a_{h}^{{\hat\ell} }(x,\theta^{\hat\ell}),\qquad {\jhat} =H_{\hat\ell}+1,\dots,{H_{\ell}+K},	\\ [0.3truecm]
	%\label{a-bar-5}
	&& \hat a_j^{{\hat\ell+1} }(x,\hat\theta^{\hat\ell+1})=\sum_{i=1, i\not=h}^{H_{{\hat\ell}}}  w_{ji}^{{\hat\ell}+1 }\;{g( \hat a_i^{{{\hat\ell}}}(x,\hat\theta^{{\hat\ell} }))}+\lambda_0 w_{jh}^{{\hat\ell}+1 }\;{g( \hat a_h^{{{\hat\ell}}}(x,\hat\theta^{{\hat\ell} }))}+\\ [0.3truecm]
	\noindent
	&&\qquad\qquad\sum_{\ihat=H_{{\hat\ell}}+1}^{{H_{\ell}+K}}  \lambda_{{\ihat-H_{\hat\ell}}} w_{j h}^{\hat\ell+1}\;{g(a_{h}^{{\hat\ell} }(x,\theta^{\hat\ell}))}+ \sigma^{{\hat\ell+1} }_{j},	\quad {j} =1,\dots,H_{\hat\ell+1}.
\end{eqnarray*}
Now, the properties of the scalars $\lambda_i$, $i=0,\ldots,K$ and   (\ref{eq28}) imply
\[
	\hat a_j^{{\hat\ell+1} }(x,\hat\theta^{\hat\ell+1})= a_j^{{\hat\ell+1} }(x,\theta^{\hat\ell+1}),\quad {j} =1,\dots,H_{\hat\ell+1},
\]
and for $\ell =\hat\ell+2,\dots,L$

\begin{eqnarray}
	\nonumber
	&&\qquad \hat a_j^{{\ell} }(x,\hat\theta^\ell)=\sum_{i=1}^{H_{\ell-1}}  w_{ji}^{{\ell} }\;{g( a_i^{{\ell} -1}(x,\theta^{{\ell} -1}))}+ \sigma^{{\ell} }_{j}	\\ [0.3truecm]
	\nonumber
	&&\qquad\qquad\qquad\qquad= a_j^{{\ell} }(x,\theta^\ell),\qquad {j} =1,\dots,H_\ell,
\end{eqnarray}	
\[
	 \hat f_r^{}(x, \hat\theta^{L})=f_r^{}(x, \theta^{L}), \qquad {r} =1,\dots,m,
\]
which hence proves (\ref{eq27}) and concludes the proof.$\hfill\Box$
\par\bigskip\noindent
\begin{proposition}\label{Prop2}
	For every  point  $\theta$, we have that, for $r=1,\ldots,m$, and  ${\ell} =1,\dots,L$
\begin{eqnarray}
	\label{eq30}
	&& \hat f_r^{}(x,\beta_{\zeta {s}}(\theta)^{L})= \hat f_r^{}(x,\gamma_{\lambda }(\theta)^{L})=f_r^{}(x, \theta^{L}), \\ [0.3truecm]
	\label{eq31}
	&&\hat a_j^{{\ell} }(x,\beta_{\zeta {s}}(\theta)^\ell)=\hat a_j^{{\ell} }(x,\gamma_{\lambda }(\theta)^\ell)= a_j^{{\ell} }(x,\theta^\ell),\quad {j} =1,\dots,H_\ell,
\end{eqnarray}
and
\begin{eqnarray}
	\label{eq32}
	&& \hat a_j^{{\hat\ell} }(x,\beta_{\zeta {s}}(\theta)^{\hat\ell})=\zeta_j ,\qquad\qquad {j} =H_{\hat\ell}+1,\dots,{H_{\ell}+K},\\ 	[0.3truecm]
\label{eq33}
	&&\hat a_j^{{\hat\ell} }(x,\gamma_{\lambda }(\theta)^{\hat\ell})	
	= a_h^{{\hat\ell} }(x,\theta^{\hat\ell}),\quad\; {j} =H_{\hat\ell}+1,\dots,{H_{\ell}+K},	
\end{eqnarray}
where $h\in [1,\dots,H_{\bar\ell}]$.
\end{proposition}
\par\medskip\noindent
{\bf Proof}. The proof follows from the one of Proposition \ref{Prop1}.$\hfill\Box$
\par\bigskip\noindent
\begin{proposition}\label{Prop3}
	Let the point  $\theta$ be such that
	\begin{eqnarray}\label{eq34}
&&\nabla_\theta R_{emp}(\theta)=0,
\end{eqnarray} 
and let the point $\hat\theta$ be given by 
\begin{eqnarray}\label{eq35}
	&&\hat\theta =\beta_{\zeta {0}}(\theta),
\end{eqnarray}  
or  
\begin{eqnarray}\label{eq36}
&&		\hat\theta =\gamma_{\lambda}(\theta).
\end{eqnarray}
Then, it results: 
	\begin{eqnarray}\label{eq37}
		&& \nabla_{\hat\theta} \hat R_{emp}(\hat\theta)=0.
	\end{eqnarray}
\end{proposition}
\par\medskip\noindent
{\bf Proof}.
The proof follows by evaluating the partial derivatives
$$\frac{\partial \hat R_{emp}(\hat\theta^L)}{ \sigma_j^\ell},\qquad \frac{\partial \hat R_{emp}(\hat\theta^L)}{w_{j,i}^\ell},$$
for $\quad {\ell} =1,\dots,L$, $\quad j=1,\ldots,\hat H_\ell$, $\quad i=1,\ldots,\hat H_{\ell-1}$ by using (\ref{eq3})-(\ref{eq10}).
\par\medskip\noindent
For the sake of brevity, we consider only the case in which $\hat\ell \in (1,L)$. The cases $\hat\ell=1$ and $\hat\ell=L$ follow from similar reasoning. 
\par\medskip\noindent
The proof is divided in the following three parts:
\begin{itemize}
	\item[a)] $\ell\in[ \hat\ell+1,L)$;
	\item[b)] $\ell=\hat\ell$;
	\item[c)] $\ell\in(1, \hat\ell-1]$.	
\end{itemize}
\par\bigskip\noindent
Part a):\quad  $\ell\in[ \hat\ell+1,L]$.
\par\bigskip\noindent
For  $j=1,\dots,\bar H_\ell$, $i=1,\ldots,H_{\ell-1}$, we can recall formulas   (\ref{eq3})-(\ref{eq10}). 
\par\medskip\noindent
Therefore, for $\tilde c=1,\ldots, \bar H_\ell$ we have:
\begin{eqnarray*}
	\label{der3}
	&&\frac{\partial \hat a^{\ell}_c(x^{p}, \hat\theta^\ell))}{\partial \hat \sigma^{{\ell}}_{j}}=\delta_{cj},\\ [0.3truecm]
	\label{der4}
	&&\frac{\partial \hat a^{\ell}_j (x^{p}, \hat\theta^\ell))}{\partial \hat w^{{\ell}}_{ji}}=	 g\biggl(\hat a_i^{\ell-1}(x^{p},\hat\theta^{\ell-1})\biggr)\delta_{cj},\\ [0.3truecm]
	&& \hbox{for}\quad {q} =\ell+1,\dots,L, \quad c=1,\ldots,H_q,	\\ [0.3truecm]
	&&\quad \frac{\partial\hat  a^{q}_c(x^{p}, \hat\theta^{q}))}{\partial \hat \sigma^{{\ell}}_{j}}=	\sum_{h=1}^{H_{q-1}}\hat w^q_{ch} g^\prime\biggl(\hat a_{h}^{q-1}(x^{p},\hat\theta^{q-1})\biggr) \frac{\partial \hat a_{h}^{q-1}(x^{p},\hat \theta^{q-1})}{\partial \hat \sigma^{{\ell}}_{j}},\\ [0.3truecm]
	&&\quad \frac{\partial \hat a^{q}_c(x^{p}, \hat\theta^q))}{\partial \hat w^{{\ell}}_{ji}}=	\sum_{h=1}^{H_{q-1}}\hat w^q_{ch} g^\prime\biggl(\hat a_{h}^{q-1}(x^{p},\hat \theta^{q-1})\biggr) \frac{\partial \hat a_{h}^{q-1}(x^{p},\hat\theta^{q-1})}{\partial \hat w^{{\ell}}_{ji}},\\ [0.3truecm]
	&&\quad \frac{\partial \hat R_{emp}(\hat\theta^L)}{\partial \hat \sigma^{{\ell}}_{j}}  = \frac{1}{P} \sum_{p=1}^{P} \sum_{r=1}^{m} { \mathcal{L}^\prime(y^{p}, \hat f_h(x^{p}, \hat \theta^L))}
	\frac{\partial \hat a^{L}_r(x^{p}, \hat \theta^L))}{\partial \hat \sigma^{{\ell}}_{j}},\\ [0.3truecm]
	&& \quad \frac{\partial \hat R_{emp}(\hat \theta^L)}{\partial \hat w^{{\ell}}_{ji}}  = \frac{1}{P} \sum_{p=1}^{P} \sum_{r=1}^{m}{ \mathcal{L}^\prime(y^{p}, \hat f_r(x^{p}, \hat\theta^L))} 
	\frac{\partial \hat a^{L}_r(x^{p}, \hat \theta^L))}{\partial \hat w^{{\ell}}_{ji}}.
\end{eqnarray*}
The definition of the maps $\beta_{\zeta 0}$ and $\gamma_\lambda$,  the equalities  (\ref{eq30}), (\ref{eq31}) and the assumption (\ref{eq34}) imply that for $\tilde c=1,\ldots, H_\ell:$
\begin{eqnarray}
	\nonumber
	&&\frac{\partial \hat a^{\ell}_c(x^{p}, \hat\theta^\ell))}{\partial \hat \sigma^{{\ell}}_{j}}=\delta_{cj}=\frac{\partial  a^{\ell}_c(x^{p}, \theta^\ell))}{\partial  \sigma^{{\ell}}_{j}},\\ [0.3truecm]
	\nonumber
	&&\frac{\partial \hat a^{\ell}_c (x^{p}, \hat\theta^\ell))}{\partial  \hat w^{{\ell}}_{ji}}=	 g\biggl( a_i^{\ell-1}(x^{p},\theta^{\ell-1})\biggr)\delta_{cj}=\frac{\partial  a^{\ell}_c (x^{p}, \theta^\ell))}{\partial  w^{{\ell}}_{ji}},\\ [0.3truecm]
	\nonumber
	&& \hbox{for}\quad {q} =\ell+1,\dots,L, \quad c=1,\ldots,H_q	\\ [0.3truecm]
	\nonumber
	&&\quad \frac{\partial\hat  a^{q}_c(x^{p}, \hat\theta^{q}))}{\partial \hat \sigma^{{\ell}}_{j}}=	\sum_{h=1}^{H_{q-1}} w^q_{ch} g^\prime\biggl( a_{h}^{q-1}(x^{p},\theta^{q-1})\biggr) \frac{\partial  a_{h}^{q-1}(x^{p}, \theta^{q-1})}{\partial  \sigma^{{\ell}}_{j}}= \frac{\partial  a^{q}_c(x^{p}, \theta^{q}))}{\partial   \sigma^{{\ell}}_{j}},\\ [0.3truecm]
	\nonumber
	&&\quad \frac{\partial \hat a^{q}_c(x^{p}, \hat\theta^q))}{\partial \hat w^{{\ell}}_{ji}}=	\sum_{h=1}^{H_{q-1}} w^q_{ch} g^\prime\biggl( a_{h}^{q-1}(x^{p}, \theta^{q-1})\biggr) \frac{\partial  a_{h}^{q-1}(x^{p},\theta^{q-1})}{\partial  w^{{\ell}}_{ji}}=\frac{\partial  a^{q}_c(x^{p}, \theta^q))}{\partial w^{{\ell}}_{ji}},\\ [0.3truecm]
	\label{eq38}
	&& \hspace*{-0.1cm} \frac{\partial \hat R_{emp}(\hat\theta^L)}{\partial \hat \sigma^{{\ell}}_{j}}  = \frac{1}{P} \sum_{p=1}^{P} \sum_{r=1}^{m} { \mathcal{L}^\prime(y^{p},  f_h(x^{p}, \theta^L))}
	\frac{\partial \hat a^{L}_r(x^{p}, \theta^L))}{\partial   \sigma^{{\ell}}_{j}}=\frac{\partial  R_{emp}(\theta^L)}{\partial   \sigma^{{\ell}}_{j}}  = 0,\\ [0.3truecm]
	\label{eq39}
	&&\hspace*{-0.1cm} \frac{\partial \hat R_{emp}(\hat \theta^L)}{\partial \hat w^{{\ell}}_{ji}}  = \frac{1}{P} \sum_{p=1}^{P} \sum_{r=1}^{m}{ \mathcal{L}^\prime(y^{p},  f_r(x^{p}, \theta^L))} 
	\frac{\partial  a^{L}_r(x^{p},  \theta^L))}{\partial  w^{{\ell}}_{ji}}=\frac{\partial  R_{emp}( \theta^L)}{\partial  w^{{\ell}}_{ji}}  =0.
\end{eqnarray}
If $\ell=\hat\ell+1$, for the components  $j=1,\dots,\hat H_{\hat\ell+1}$, and $i=H_{\hat\ell}+1,\ldots,{H_{\hat\ell}+K}$ we have that:
\par\medskip\noindent 
when $\hat\theta =\beta_{\zeta {0}}(\theta)$, the property  (\ref{eq19}) implies
\par\medskip\noindent 
\begin{eqnarray}
	\nonumber
	&& \hbox{for }\quad  \tilde c=1,\ldots,H_{\hat\ell+1}	\\ [0.3truecm]
	\nonumber
	&&\frac{\partial \hat a^{\hat\ell+1}_{\tilde c} (x^{p}, \hat\theta^{\hat\ell+1}))}{\partial \hat w^{{\hat\ell+1}}_{ji}}=	 g\bigl(\zeta_i\bigr)\delta_{\tilde c j}=g\bigl(\zeta_i\bigr)\frac{\partial \hat a^{\hat\ell+2}_{\tilde c }(x^{p}, \hat\theta^{\hat\ell+2})}{\partial \hat \sigma^{\hat\ell+1}_{j}},\\ [0.3truecm]
	\nonumber
%	&& \hbox{for }\quad  \tilde c=1,\ldots,H_{\hat\ell+2}	\\ [0.3truecm]
%	\nonumber
%	&&\quad\frac{\partial \hat a^{\hat\ell+2}_{\tilde c }(x^{p}, \hat\theta^{\hat\ell+2})}{\partial \hat w^{\hat\ell+1}_{ji}}=g\bigl(\zeta_i\bigr)\sum_{h=1}^{H_{{\hat\ell+1}}}\hat w^{\hat\ell+2}_{{\tilde c }h} g^\prime\biggl(\hat a_{h}^{{\hat\ell+1}}(x^{p},\hat \theta^{{\hat\ell+1}})\biggr)=g\bigl(\zeta_i\bigr) \frac{\partial \hat a^{\hat\ell+2}_{\tilde c }(x^{p}, \hat\theta^{\hat\ell+2})}{\partial \hat \sigma^{\hat\ell+1}_{j}}\\ [0.3truecm]
%	\nonumber
	&& \hbox{for }\quad {q} =\hat\ell+2,\dots,L, \quad c=1,\ldots,H_q	\\ [0.3truecm]
	\nonumber
	&&\quad \frac{\partial \hat a^{q}_c(x^{p}, \hat\theta^q))}{\partial \hat w^{{\hat\ell+1}}_{ji}}=g\bigl(\zeta_i\bigr)	\sum_{h=1}^{H_{q-1}}\hat w^q_{ch} g^\prime\biggl(\hat a_{h}^{q-1}(x^{p},\hat \theta^{q-1})\biggr) \frac{\partial \hat a_{h}^{q-1}(x^{p},\hat\theta^{q-1})}{\partial \hat \sigma^{\hat\ell+1}_{j}}=g\bigl(\zeta_i\bigr)\frac{\partial \hat a^{q}_c(x^{p}, \hat\theta^q))}{\partial \hat \sigma^{\hat\ell+1}_{j}},\\ [0.5truecm]
	\label{eq40}
	&&\frac{\partial \hat R_{emp}(\hat \theta^L)}{\partial \hat w^{{\hat\ell+1}}_{ji}}  = g\bigl(\zeta_i\bigr)\frac{1}{P} \sum_{p=1}^{P} \sum_{r=1}^{m}{ \mathcal{L}^\prime(y^{p}, \hat f_r(x^{p}, \hat\theta^L))} 
	\frac{\partial \hat a^{L}_r(x^{p}, \hat \theta^L))}{\partial \hat \sigma^{\hat\ell+1}_{j}}\\ \nonumber
&&\qquad\qquad\quad\;\;= g\bigl(\zeta_i\bigr)	\frac{\partial \hat R_{emp}(\hat \theta^L)}{\partial \hat \sigma^{{\hat\ell+1}}_{j}}=0,
\end{eqnarray}
where the last equality follows from (\ref{eq38}).
\par\medskip\noindent
Instead, when $\hat\theta =\gamma_{\lambda}(\theta)$ the property  (\ref{eq21}) ensures that there exist an index 
$$ z\in [1,\ldots,H_{\hat\ell}],$$ 
such that  
\begin{eqnarray}
	\nonumber
	&& \hbox{for }\quad  \tilde c=1,\ldots,H_{\hat\ell+1}	\\ [0.3truecm]
	\nonumber
	&&\quad \frac{\partial \hat a^{\hat\ell+1}_{\tilde c} (x^{p}, \hat\theta^{\hat\ell+1}))}{\partial \hat w^{{\hat\ell+1}}_{ji}}=	 g\bigl(\hat a_z^{{\hat\ell} }(x^p,\hat\theta^{\hat\ell})\bigr)\delta_{\tilde c j}=\frac{\partial \hat a^{\hat\ell+1}_{\tilde c} (x^{p}, \hat\theta^{\hat\ell+1}))}{\partial \hat w^{{\hat\ell+1}}_{jz}},\\ [0.3truecm]
	\nonumber
	&& \hbox{for }\quad {q} =\hat\ell+2,\dots,L, \quad c=1,\ldots,H_q	\\ [0.3truecm]
	\nonumber
	&&\quad \frac{\partial \hat a^{q}_c(x^{p}, \hat\theta^q))}{\partial \hat w^{{\ell}}_{ji}}=	\sum_{h=1}^{H_{q-1}}\hat w^q_{ch} g^\prime\biggl(\hat a_{h}^{q-1}(x^{p},\hat \theta^{q-1})\biggr) \frac{\partial \hat a_{h}^{q-1}(x^{p},\hat\theta^{q-1})}{\partial \hat w^{\hat\ell+1}_{jz}}=\frac{\partial \hat a^{q}_c(x^{p}, \hat\theta^q))}{\partial \hat w^{{\ell}}_{jz}},\\ [0.5truecm]
	\label{eq41}
	&&\hspace*{-0.1cm} \frac{\partial \hat R_{emp}(\hat \theta^L)}{\partial \hat w^{{\ell}}_{ji}}  = \frac{1}{P} \sum_{p=1}^{P} \sum_{r=1}^{m}{ \mathcal{L}^\prime(y^{p}, \hat f_r(x^{p}, \hat\theta^L))} 
	\frac{\partial \hat a^{L}_r(x^{p}, \hat \theta^L))}{\partial \hat w^{{\ell}}_{jz}}=\frac{\partial \hat R_{emp}(\hat \theta^L)}{\partial \hat w^{{\ell}}_{jz}}=0,
\end{eqnarray}
where (\ref{eq39}) gives the last equality.
\par\bigskip\noindent
Part b):\quad  ${\ell} =\hat\ell$.
\par\medskip\noindent
For every   $j=1,\ldots,H_{\hat\ell},$ and $ i=1,\ldots,\hat H_{\hat\ell-1}$ the expressions  (\ref{eq3})-(\ref{eq10}) yield:
\par\medskip\noindent
\begin{eqnarray}
	\nonumber 
	&& \hbox{for }\quad  \tilde c=1,\ldots,\hat H_{\hat\ell}	\\ [0.3truecm]
	\nonumber
	&&\quad \hbox{for }\quad  j=1,\ldots,H_{\hat\ell},\quad   i=1,\ldots,\hat H_{\hat\ell-1}	\\ [0.3truecm]
	\nonumber		
	&&\qquad\frac{\partial \hat a^{\hat\ell}_{\tilde c}(x^{p}, \hat\theta^{\hat\ell}))}{\partial \hat \sigma^{{\hat\ell}}_{j}}=\delta_{\tilde c j}=\frac{\partial  a^{\hat\ell}_{\tilde c}(x^{p}, \theta^{\hat\ell}))}{\partial  \sigma^{{\hat\ell}}_{{j}}},\\ [0.3truecm]
	\nonumber
	&&\qquad\frac{\partial \hat a^{\hat\ell}_{\tilde c} (x^{p}, \hat\theta^{\hat\ell}))}{\partial \hat w^{{\hat\ell}}_{ji}}=	 g\biggl(\hat a_i^{{\hat\ell}-1}(x^{p},\hat\theta^{{\hat\ell}-1})\biggr)\delta_{\tilde c j}=\frac{\partial  a^{\hat\ell}_{\tilde c} (x^{p}, \theta^{\hat\ell}))}{\partial  w^{{\hat\ell}}_{ji}}, 
		\\ [0.3truecm]
	\nonumber
	&&\quad \hbox{for }\quad  j=H_{\hat\ell}+1,\ldots,{H_{\ell}+K},\quad   i=1,\ldots,\hat H_{\hat\ell-1}\\ [0.3truecm]
	\nonumber	
	&&\qquad\frac{\partial \hat a^{\hat\ell}_{{\tilde c}}(x^{p}, \hat\theta^{\hat\ell}))}{\partial \hat \sigma^{{\hat\ell}}_{\jhat}}=\delta_{\tilde c \jhat},\\ [0.3truecm]
	\nonumber
	&&\qquad\frac{\partial \hat a^{\hat\ell}_{{\tilde c}} (x^{p}, \hat\theta^{\hat\ell}))}{\partial \hat w^{{\hat\ell}}_{{\jhat}i}}=	 g\biggl(\hat a_i^{{\hat\ell}-1}(x^{p},\hat\theta^{{\hat\ell}-1})\biggr)\delta_{\tilde c \jhat}.
	\end{eqnarray}
\par\medskip\noindent
When $\hat\theta =\beta_{\zeta {0}}(\theta)$, the properties  (\ref{eq20}) and (\ref{eq31})  yield:
%\par\medskip\noindent 
%Then for all $\tilde c=1,\ldots,\hat H_{{\hat\ell}+1}$, $j=1,\ldots,H_{\hat\ell}$ and $ i=1,\ldots,\hat H_{\hat\ell-1}$ we have:
	\begin{eqnarray}
	\nonumber 
	&& \hbox{for }\quad  \tilde c=1,\ldots,\hat H_{{\hat\ell}+1}	\\ [0.3truecm]
	\nonumber
	&&\quad \hbox{for }\quad  j=1,\ldots,H_{\hat\ell},\quad   i=1,\ldots,\hat H_{\hat\ell-1}	\\ [0.3truecm]
	\nonumber
	&&\qquad\frac{\partial \hat a^{{\hat\ell}+1}_{\tilde c}(x^{p}, \hat\theta^{{\hat\ell}+1}))}{\partial \hat \sigma^{{\hat\ell}}_{j}}=	
		 w^{{\hat\ell}+1}_{{\tilde c}j} g^\prime\biggl( a_{j}^{\hat\ell}(x^{p},\theta^{{\hat\ell}})\biggr)	 
	\frac{\partial  a_{j}^{{\hat\ell}}(x^{p},\theta^{{\hat\ell}})}{\partial   \sigma^{{\hat\ell}}_{j}}=\frac{\partial  a^{{\hat\ell}+1}_{\tilde c}(x^{p}, \theta^{{\hat\ell}+1}))}{\partial  \sigma^{{\hat\ell}}_{j}},
%	\nonumber
%&&	\qquad\qquad\qquad\qquad\qquad=\frac{\partial  a^{{\hat\ell}+1}_{\tilde c}(x^{p}, \theta^{{\hat\ell}+1}))}{\partial  \sigma^{{\hat\ell}}_{j}}
	\\ [0.3truecm]
	\nonumber
	&&\qquad\frac{\partial \hat a^{{\hat\ell}+1}_{\tilde c}(x^{p}, \hat\theta^{{\hat\ell}+1}))}{\partial \hat w^{\hat\ell}_{ji}}=	
w^{{\hat\ell}+1}_{{\tilde c}j} g^\prime\biggl( a_{j}^{\hat\ell}(x^{p},\theta^{{\hat\ell}})\biggr)	 
\frac{\partial  a_{j}^{{\hat\ell}}(x^{p},\theta^{{\hat\ell}})}{\partial   w^{\hat\ell}_{ji}}=\frac{\partial  a^{{\hat\ell}+1}_{\tilde c}(x^{p}, \theta^{{\hat\ell}+1}))}{\hat w^{\hat\ell}_{ji}},\\
[0.3truecm]
\nonumber
&&\quad \hbox{for }\quad  \jhat=H_{\hat\ell}+1,\ldots,{H_{\ell}+K},\quad   i=1,\ldots,\hat H_{\hat\ell-1}	\\ [0.3truecm]
\nonumber
&&\qquad\frac{\partial \hat a^{{\hat\ell}+1}_{\tilde c}(x^{p}, \hat\theta^{{\hat\ell}+1}))}{\partial \hat \sigma^{{\hat\ell}}_{\jhat}}=	
0\; g^\prime\bigl(\zeta_{\jhat}\bigr)	 
\frac{\partial  a_{j}^{{\hat\ell}}(x^{p},\theta^{{\hat\ell}})}{\partial   \sigma^{{\hat\ell}}_{\jhat}}=0,
\\ [0.3truecm]
\nonumber
&&\qquad\frac{\partial \hat a^{{\hat\ell}+1}_{\tilde c}(x^{p}, \hat\theta^{{\hat\ell}+1}))}{\partial \hat w^{\hat\ell}_{\jhat i}}=	
0\; g^\prime\bigl( \zeta_{\jhat}\bigr)	 
\frac{\partial  a_{j}^{{\hat\ell}}(x^{p},\theta^{{\hat\ell}})}{\partial   w^{\hat\ell}_{\jhat i}}=0.\\
[0.3truecm]
\nonumber	
\end{eqnarray}
\noindent Therefore, we have:
\begin{eqnarray}
\nonumber
&& \hbox{for } {q} ={\hat\ell}+2,\dots,L, \quad c=1,\ldots,\hat H_q	\\ [0.3truecm]
\nonumber
&&\quad \hbox{for }\quad  j=1,\ldots,H_{\hat\ell},\quad   i=1,\ldots,\hat H_{\hat\ell-1}	\\ [0.3truecm]
\nonumber
&&\qquad \frac{\partial \hat a^{q}_c(x^{p}, \hat\theta^{q}))}{\partial \hat \sigma^{{\ell}}_{j}}=	\sum_{h=1}^{H_{q-1}}w^q_{ch} g^\prime\biggl(a_{h}^{q-1}(x^{p},\theta^{q-1})\biggr) \frac{\partial a_{h}^{q-1}(x^{p},\theta^{q-1})}{\partial  \sigma^{{\ell}}_{j}},\\ [0.3truecm]
\nonumber
&&\qquad \frac{\partial \hat a^{q}_c(x^{p}, \hat \theta^q))}{\partial \hat w^{{\ell}}_{ji}}=	\sum_{h=1}^{H_{q-1}}w^q_{ch} g^\prime\biggl(a_{h}^{q-1}(x^{p},\theta^{q-1})\biggr) \frac{\partial a_{h}^{q-1}(x^{p},\theta^{q-1})}{\partial w^{{\ell}}_{ji}},\\ [0.3truecm]
	\nonumber
	&&\quad\frac{\partial \hat R_{emp}(\hat \theta^L)}{\partial \hat \sigma^{{\ell}}_{j}}  = \frac{1}{P} \sum_{p=1}^{P} \sum_{r=1}^{m} { \mathcal{L}^\prime(y^{p}, f_h(x^{p}, \theta^L))}
	\frac{\partial a^{L}_r(x^{p}, \theta^L))}{\partial  \sigma^{{\ell}}_{j}}\\
	\label{eq42}
	&&\qquad\qquad\qquad\quad=\frac{\partial  R_{emp}( \theta^L)}{\partial  \sigma^{{\ell}}_{j}}  = 0,\\ [0.3truecm]
	\nonumber
	&&\quad\frac{\partial\hat  R_{emp}(\hat \theta^L)}{\partial \hat w^{{\ell}}_{ji}}  = \frac{1}{P} \sum_{p=1}^{P} \sum_{r=1}^{m}{ \mathcal{L}^\prime(y^{p}, f_r(x^{p}, \theta^L))} 
	\frac{\partial a^{L}_r(x^{p}, \theta^L))}{\partial w^{{\ell}}_{ji}}\\
	\label{eq43}
	&&\qquad\qquad\qquad\quad=\frac{\partial  R_{emp}( \theta^L)}{\partial w^{{\ell}}_{ji}}  = 0,
	\\ [0.3truecm]
\nonumber
&&\quad \hbox{for }\quad  \jhat =H_{\hat\ell}+1,\ldots,{H_{\ell}+K},\quad   i=1,\ldots,\hat H_{\hat\ell-1}	\\ [0.3truecm]
\nonumber
&&\qquad \frac{\partial \hat a^{q}_c(x^{p}, \hat\theta^{q}))}{\partial \hat \sigma^{{\ell}}_{\jhat}}=	\sum_{h=1}^{H_{q-1}}w^q_{ch} g^\prime\biggl(a_{h}^{q-1}(x^{p},\theta^{q-1})\biggr) 0=0,\\ [0.3truecm]
\nonumber
&&\qquad \frac{\partial \hat a^{q}_c(x^{p}, \hat \theta^q))}{\partial \hat w^{{\ell}}_{{\jhat} i} }=	\sum_{h=1}^{H_{q-1}}w^q_{ch} g^\prime\biggl(a_{h}^{q-1}(x^{p},\theta^{q-1})\biggr)0=0,\\ [0.3truecm]
\label{eq44}
&&\quad\frac{\partial \hat R_{emp}(\hat \theta^L)}{\partial \hat \sigma^{{\ell}}_{\jhat}}  = \frac{1}{P} \sum_{p=1}^{P} \sum_{r=1}^{m} { \mathcal{L}^\prime(y^{p}, f_h(x^{p}, \theta^L))}
0 = 0,\\ [0.3truecm]
\label{eq45}
&&\quad\frac{\partial\hat  R_{emp}(\hat \theta^L)}{\partial \hat w^{{\ell}}_{\jhat i}}  = \frac{1}{P} \sum_{p=1}^{P} \sum_{r=1}^{m}{ \mathcal{L}^\prime(y^{p}, f_r(x^{p}, \theta^L))} 
0= 0.	
\end{eqnarray}
Instead, when  $\hat\theta =\gamma_{\lambda}(\theta)$   there exist an index 
$$ z\in [1,\ldots,H_{\hat\ell}]$$ 
such that property  (\ref{eq21})-(\ref{eq24}) hold. Thus, we have:
	\begin{eqnarray}
	\nonumber 
	&& \hbox{for }\quad  \tilde c=1,\ldots,\hat H_{{\hat\ell}+1}	\\ [0.3truecm]
	\nonumber
	&&\quad \hbox{for }\quad  j=1,\ldots,H_{\hat\ell},\quad j\not=z, \quad  i=1,\ldots,\hat H_{\hat\ell-1}	\\ [0.3truecm]
	\nonumber
	&&\qquad\frac{\partial \hat a^{{\hat\ell}+1}_{\tilde c}(x^{p}, \hat\theta^{{\hat\ell}+1}))}{\partial \hat \sigma^{{\hat\ell}}_{j}}=	
	w^{{\hat\ell}+1}_{{\tilde c}j} g^\prime\biggl( a_{j}^{\hat\ell}(x^{p},\theta^{{\hat\ell}})\biggr)	 
	\frac{\partial  a_{j}^{{\hat\ell}}(x^{p},\theta^{{\hat\ell}})}{\partial   \sigma^{{\hat\ell}}_{j}}=\frac{\partial  a^{{\hat\ell}+1}_{\tilde c}(x^{p}, \theta^{{\hat\ell}+1}))}{\partial  \sigma^{{\hat\ell}}_{j}},
	%	\nonumber
	%&&	\qquad\qquad\qquad\qquad\qquad=\frac{\partial  a^{{\hat\ell}+1}_{\tilde c}(x^{p}, \theta^{{\hat\ell}+1}))}{\partial  \sigma^{{\hat\ell}}_{j}}
	\\ [0.3truecm]
	\nonumber
	&&\qquad\frac{\partial \hat a^{{\hat\ell}+1}_{\tilde c}(x^{p}, \hat\theta^{{\hat\ell}+1}))}{\partial \hat w^{\hat\ell}_{ji}}=	
	w^{{\hat\ell}+1}_{{\tilde c}j} g^\prime\biggl( a_{j}^{\hat\ell}(x^{p},\theta^{{\hat\ell}})\biggr)	 
	\frac{\partial  a_{j}^{{\hat\ell}}(x^{p},\theta^{{\hat\ell}})}{\partial   w^{\hat\ell}_{ji}}=\frac{\partial  a^{{\hat\ell}+1}_{\tilde c}(x^{p}, \theta^{{\hat\ell}+1}))}{\hat w^{\hat\ell}_{ji}},\\
	[0.3truecm]
	\nonumber
	&&\quad \hbox{for }\quad  j=z,\quad   i=1,\ldots,\hat H_{\hat\ell-1}	\\ [0.3truecm]
	\nonumber
	&&\qquad\frac{\partial \hat a^{{\hat\ell}+1}_{\tilde c}(x^{p}, \hat\theta^{{\hat\ell}+1}))}{\partial \hat \sigma^{{\hat\ell}}_{z}}=	\lambda_0
	w^{{\hat\ell}+1}_{{\tilde c}z} g^\prime\biggl( a_{z}^{\hat\ell}(x^{p},\theta^{{\hat\ell}})\biggr)	 
	\frac{\partial  a_{z}^{{\hat\ell}}(x^{p},\theta^{{\hat\ell}})}{\partial   \sigma^{{\hat\ell}}_{z}}=\lambda_0\frac{\partial  a^{{\hat\ell}+1}_{\tilde c}(x^{p}, \theta^{{\hat\ell}+1}))}{\partial  \sigma^{{\hat\ell}}_{z}},
	%	\nonumber
	%&&	\qquad\qquad\qquad\qquad\qquad=\frac{\partial  a^{{\hat\ell}+1}_{\tilde c}(x^{p}, \theta^{{\hat\ell}+1}))}{\partial  \sigma^{{\hat\ell}}_{j}}
	\\ [0.3truecm]
	\nonumber
	&&\qquad\frac{\partial \hat a^{{\hat\ell}+1}_{\tilde c}(x^{p}, \hat\theta^{{\hat\ell}+1}))}{\partial \hat w^{\hat\ell}_{zi}}=\lambda_0
	w^{{\hat\ell}+1}_{{\tilde c}z} g^\prime\biggl( a_{z}^{\hat\ell}(x^{p},\theta^{{\hat\ell}})\biggr)	 
	\frac{\partial  a_{z}^{{\hat\ell}}(x^{p},\theta^{{\hat\ell}})}{\partial   w^{\hat\ell}_{zi}}=\lambda_0\frac{\partial  a^{{\hat\ell}+1}_{\tilde c}(x^{p}, \theta^{{\hat\ell}+1}))}{\hat w^{\hat\ell}_{zi}},\\
	[0.3truecm]
	\nonumber
	&&\quad \hbox{for }\quad  \jhat=H_{\hat\ell}+1,\ldots,{H_{\ell}+K},\quad   i=1,\ldots,\hat H_{\hat\ell-1}	\\ [0.3truecm]
	\nonumber
	&&\qquad\frac{\partial \hat a^{{\hat\ell}+1}_{\tilde c}(x^{p}, \hat\theta^{{\hat\ell}+1}))}{\partial \hat \sigma^{{\hat\ell}}_{\jhat}}=\lambda_{\jhat-H_{\hat\ell}}
	w^{{\hat\ell}+1}_{{\tilde c}z} g^\prime\biggl( a_{z}^{\hat\ell}(x^{p},\theta^{{\hat\ell}})\biggr)	 
	\frac{\partial  a_{z}^{{\hat\ell}}(x^{p},\theta^{{\hat\ell}})}{\partial   \sigma^{{\hat\ell}}_{z}}=\lambda_{\jhat-H_{\hat\ell}}\frac{\partial  a^{{\hat\ell}+1}_{\tilde c}(x^{p}, \theta^{{\hat\ell}+1}))}{\partial  \sigma^{{\hat\ell}}_{z}},\\
	[0.3truecm]
	&&\qquad\frac{\partial \hat a^{{\hat\ell}+1}_{\tilde c}(x^{p}, \hat\theta^{{\hat\ell}+1}))}{\partial \hat w^{\hat\ell}_{\jhat i}}=\lambda_{\jhat-H_{\hat\ell}}	
	w^{{\hat\ell}+1}_{{\tilde c}z} g^\prime\biggl( a_{z}^{\hat\ell}(x^{p},\theta^{{\hat\ell}})\biggr)	 
	\frac{\partial  a_{z}^{{\hat\ell}}(x^{p},\theta^{{\hat\ell}})}{\partial   w^{\hat\ell}_{zi}}=\lambda_{\jhat-H_{\hat\ell}}\frac{\partial  a^{{\hat\ell}+1}_{\tilde c}(x^{p}, \theta^{{\hat\ell}+1}))}{\hat w^{\hat\ell}_{zi}}.
	\nonumber	
\end{eqnarray}
Now, by recalling that $\sum_{{\tilde j}=0}^{K}\lambda_{\tilde j}=1$ and by repeating the same reasoning for the case $\hat\theta =\beta_{\zeta 0}(\theta)$, we obtain again (\ref{eq44}) and (\ref{eq45}).
\par\bigskip\noindent
Part c):\quad  $\ell\in(1,\ldots \hat\ell-1]$.
\par\bigskip\noindent
For  $\ell\in(1,\ldots \hat\ell-1]$,\quad $j=1,\dots,\bar H_\ell$,\quad $i=1,\ldots,\hat H_{\ell-1}$ we can use again formulas   (\ref{eq3})-(\ref{eq10}), the definition of the maps $\beta_{\zeta 0}$ and $\gamma_\lambda$,  the equalities  (\ref{eq30}) and (\ref{eq31}). 
\par\medskip\noindent
In particular, we can write:
\par\medskip\noindent
for\quad $\tilde c=1,\ldots\hat H_\ell$\\
\begin{eqnarray}
	%	\nonumber
	%	&& \hbox{for}\quad \tilde c=1,\ldots\hat H_\ell:\\
	%	[0.3truecm]
	\nonumber
	&&\frac{\partial \hat a^{\ell}_c(x^{p}, \hat\theta^\ell))}{\partial \hat \sigma^{{\ell}}_{j}}=\delta_{cj}=\frac{\partial  a^{\ell}_c(x^{p}, \theta^\ell))}{\partial  \sigma^{{\ell}}_{j}},\\ [0.3truecm]
	\nonumber
	&&\frac{\partial \hat a^{\ell}_c (x^{p}, \hat\theta^\ell))}{\partial  \hat w^{{\ell}}_{ji}}=	 g\biggl( a_i^{\ell-1}(x^{p},\theta^{\ell-1})\biggr)\delta_{cj}=\frac{\partial  a^{\ell}_c (x^{p}, \theta^\ell))}{\partial  w^{{\ell}}_{ji}},
\end{eqnarray}
for\quad \quad ${q} =\ell+1,\dots,L$
\begin{eqnarray}
	\nonumber
	&&\hbox{if }\quad {q} \not= \hat\ell\quad \hbox{and}\quad {q} \not= \hat\ell +1, \quad \hbox{for }\quad  c=1,\ldots,\hat H_q		\\
	[0.3truecm]
	\nonumber
	&&\quad \frac{\partial\hat  a^{q}_c(x^{p}, \hat\theta^{q}))}{\partial \hat \sigma^{{\ell}}_{j}}=	\sum_{h=1}^{H_{q-1}} w^q_{ch} g^\prime\biggl( a_{h}^{q-1}(x^{p},\theta^{q-1})\biggr) \frac{\partial  a_{h}^{q-1}(x^{p}, \theta^{q-1})}{\partial  \sigma^{{\ell}}_{j}}= \frac{\partial  a^{q}_c(x^{p}, \theta^{q}))}{\partial   \sigma^{{\ell}}_{j}},\\ [0.3truecm]
	\nonumber
	&&\quad \frac{\partial \hat a^{q}_c(x^{p}, \hat\theta^q))}{\partial \hat w^{{\ell}}_{ji}}=	\sum_{h=1}^{H_{q-1}} w^q_{ch} g^\prime\biggl( a_{h}^{q-1}(x^{p}, \theta^{q-1})\biggr) \frac{\partial  a_{h}^{q-1}(x^{p},\theta^{q-1})}{\partial  w^{{\ell}}_{ji}}=\frac{\partial  a^{q}_c(x^{p}, \theta^q))}{\partial w^{{\ell}}_{ji}},\\ [0.3truecm]
		\nonumber
	&&\hbox{if }\quad {q} = \hat\ell,\quad  \quad \hbox{for }\quad  c=1,\ldots, H_{\hat\ell}		\\
	[0.3truecm]
	\nonumber
	&&\quad \frac{\partial\hat  a^{{\hat\ell}}_c(x^{p}, \hat\theta^{{\hat\ell}}))}{\partial \hat \sigma^{{\ell}}_{j}}=	\sum_{h=1}^{H_{{\hat\ell}-1}} w^{\hat\ell}_{ch} g^\prime\biggl( a_{h}^{{\hat\ell}-1}(x^{p},\theta^{{\hat\ell}-1})\biggr) \frac{\partial  a_{h}^{{\hat\ell}-1}(x^{p}, \theta^{{\hat\ell}-1})}{\partial  \sigma^{{\ell}}_{j}}= \frac{\partial  a^{{\hat\ell}}_c(x^{p}, \theta^{{\hat\ell}}))}{\partial   \sigma^{{\ell}}_{j}},\\ [0.3truecm]
	\nonumber
	&&\quad \frac{\partial \hat a^{{\hat\ell}}_c(x^{p}, \hat\theta^{\hat\ell}))}{\partial \hat w^{{\ell}}_{ji}}=	\sum_{h=1}^{H_{{\hat\ell}-1}} w^{\hat\ell}_{ch} g^\prime\biggl( a_{h}^{{\hat\ell}-1}(x^{p}, \theta^{{\hat\ell}-1})\biggr) \frac{\partial  a_{h}^{q-1}(x^{p},\theta^{q-1})}{\partial  w^{{\ell}}_{ji}}=\frac{\partial  a^{{\hat\ell}}_c(x^{p}, \theta^{\hat\ell}))}{\partial w^{{\ell}}_{ji}},
\end{eqnarray}
\begin{eqnarray}	
%	\nonumber
%	&& \hbox{for }\quad {q} =\ell+1,\dots,L:\\
%	[0.3truecm]
	\nonumber
	&&\hbox{if }\quad {q} = \hat\ell\quad \hbox{and}\quad \hat \theta=\beta_{\zeta 0}(\theta), \quad \hbox{for }\quad  {\hat c}=H_{\hat\ell}+1,\ldots, {H_{\hat\ell}+K}		\\
	[0.3truecm]
	\nonumber
	&&\quad \frac{\partial\hat  a^{{\hat\ell}}_{\hat c}(x^{p}, \hat\theta^{{\hat\ell}}))}{\partial \hat \sigma^{{\ell}}_{j}}=	\sum_{h=1}^{H_{{\hat\ell}-1}}\hat  w^{\hat\ell}_{{\hat c}h} g^\prime\biggl( a_{h}^{{\hat\ell}-1}(x^{p},\theta^{{\hat\ell}-1})\biggr) \frac{\partial  a_{h}^{{\hat\ell}-1}(x^{p}, \theta^{{\hat\ell}-1})}{\partial  \sigma^{{\ell}}_{j}},\\ [0.3truecm]
	\nonumber
	&&\quad \frac{\partial \hat a^{{\hat\ell}}_{\hat c}(x^{p}, \hat\theta^{\hat\ell}))}{\partial \hat w^{{\ell}}_{ji}}=	\sum_{h=1}^{H_{{\hat\ell}-1}} \hat w^{\hat\ell}_{{\hat c}h} g^\prime\biggl( a_{h}^{{\hat\ell}-1}(x^{p}, \theta^{{\hat\ell}-1})\biggr) \frac{\partial  a_{h}^{\hat\ell-1}(x^{p},\theta^{\hat\ell-1})}{\partial  w^{{\ell}}_{ji}},\qquad\qquad\quad
\end{eqnarray}

\begin{eqnarray}
\nonumber	&&\hbox{if }\quad {q} = {\hat\ell+1}\quad \hbox{and}\quad \hat \theta=\beta_{\zeta 0}(\theta), \quad \hbox{for }\quad  {\tilde c}=1,\ldots, H_{\hat\ell+1}		\\
	[0.3truecm]
	\nonumber
&&\frac{\partial \hat a^{{\hat\ell}+1}_{\tilde c}(x^{p}, \hat\theta^{{\hat\ell}+1}))}{\partial \hat \sigma^{{\hat\ell}}_{j}}=	
\sum_{h=1}^{H_{\hat\ell}} w^{{\hat\ell}+1}_{{\tilde c}h} g^\prime\biggl( a_{h}^{\hat\ell}(x^{p},\theta^{{\hat\ell}})\biggr)	 
\frac{\partial  a_{h}^{{\hat\ell}}(x^{p},\theta^{{\hat\ell}})}{\partial   \sigma^{{\hat\ell}}_{j}}\\
\nonumber 
&&	\qquad\qquad\qquad\qquad\qquad+\sum_{\hhat=H_{\hat\ell}+1}^{{H_{\hat\ell}+K}}0\; g^\prime\biggl(\hat a_{\hhat}^{\hat\ell}(x^{p},\hat\theta^{{\hat\ell}})\biggr)	 
\frac{\partial \hat a_{\hhat}^{{\hat\ell}}(x^{p},\hat\theta^{{\hat\ell}})}{\partial \hat \sigma^{{\ell}}_{j}}
\\ [0.3truecm]
\nonumber
&&	\qquad\qquad\qquad\qquad\qquad=\frac{\partial  a^{{\hat\ell}+1}_{\tilde c}(x^{p}, \theta^{{\hat\ell}+1}))}{\partial  \sigma^{{\hat\ell}}_{j}},
\\ [0.3truecm]
\nonumber
&&\frac{\partial \hat a^{{\hat\ell}+1}_{\tilde c}(x^{p}, \hat\theta^{{\hat\ell}+1}))}{\partial \hat w^{\hat\ell}_{j i}}=	
\sum_{h=1}^{H_{\hat\ell}} w^{{\hat\ell}+1}_{{\tilde c}h} g^\prime\biggl( a_{h}^{\hat\ell}(x^{p},\theta^{{\hat\ell}})\biggr)	 
\frac{\partial  a_{h}^{{\hat\ell}}(x^{p},\theta^{{\hat\ell}})}{\partial  w^{\hat\ell}_{j i}}\\
\nonumber
&&	\qquad\qquad\qquad\qquad\qquad+\sum_{\hhat=H_{\hat\ell}+1}^{{H_{\hat\ell}+K}}\hat 0\; g^\prime\biggl(j a_{\hhat}^{\hat\ell}(x^{p},\hat\theta^{{\hat\ell}})\biggr)	 
\frac{\partial \hat a_{\hhat}^{{\hat\ell}}(x^{p},\hat\theta^{{\hat\ell}})}{\partial \hat w^{\hat\ell}_{j i}}\\ [0.3truecm]
&&	\qquad\qquad\qquad\qquad\qquad=\frac{\partial  a^{{\hat\ell}+1}_{\tilde c}(x^{p}, \theta^{{\hat\ell}+1}))}{\partial w^{{\hat\ell}}_{ji}},
\nonumber
\end{eqnarray}

\begin{eqnarray}	
	%	\nonumber
	%	&& \hbox{for }\quad {q} =\ell+1,\dots,L:\\
	%	[0.3truecm]
	\nonumber
	&&\hbox{if }\quad {q} = \hat\ell\quad \hbox{and}\ \hat \theta=\gamma_{\lambda}(\theta),\  \hbox{for }\ z\in [1,\ldots,H_\ell]\ \hbox{and}  \ \hbox{for }\  {\hat c}=H_{\hat\ell}+1,\ldots,H_{\hat\ell}+K		\\
	[0.3truecm]
	\nonumber
	&&\quad \frac{\partial\hat  a^{{\hat\ell}}_{\hat c}(x^{p}, \hat\theta^{{\hat\ell}}))}{\partial \hat \sigma^{{\ell}}_{j}}=	\sum_{h=1}^{H_{{\hat\ell}-1}}  w^{\hat\ell}_{{z}h} g^\prime\biggl( a_{h}^{{\hat\ell}-1}(x^{p},\theta^{{\hat\ell}-1})\biggr) \frac{\partial  a_{h}^{{\hat\ell}-1}(x^{p}, \theta^{{\hat\ell}-1})}{\partial  \sigma^{{\ell}}_{j}}=\frac{\partial  a^{{\hat\ell}}_{\hat c}(x^{p}, \theta^{{\hat\ell}}))}{\partial   \sigma^{{\ell}}_{j}},\\ [0.3truecm]
	\nonumber
	&&\quad \frac{\partial \hat a^{{\hat\ell}}_{\hat c}(x^{p}, \hat\theta^{\hat\ell}))}{\partial \hat w^{{\ell}}_{ji}}=	\sum_{h=1}^{H_{{\hat\ell}-1}}  w^{\hat\ell}_{{z}h} g^\prime\biggl( a_{h}^{{\hat\ell}-1}(x^{p}, \theta^{{\hat\ell}-1})\biggr) \frac{\partial  a_{h}^{{\hat\ell}-1}(x^{p},\theta^{{\hat\ell}-1})}{\partial  w^{{\ell}}_{ji}}=\frac{\partial  a^{{\hat\ell}}_{\hat c}(x^{p}, \theta^{{\hat\ell}}))}{\partial  w^{{\ell}}_{ji}},
\end{eqnarray}

\begin{eqnarray}
	\nonumber	&&\hbox{if }\quad {q} = {\hat\ell+1}\quad \hbox{and}\quad \hat \theta=\gamma_{\lambda}(\theta),\quad  \hbox{for }\ z\in [1,\ldots,H_\ell]\quad \hbox{and}  \quad \hbox{for }\quad  {\tilde c}=1,\ldots, H_{\hat\ell+1}		\\
	[0.3truecm]
	\nonumber
	&&\frac{\partial \hat a^{{\hat\ell}+1}_{\tilde c}(x^{p}, \hat\theta^{{\hat\ell}+1}))}{\partial \hat \sigma^{{\hat\ell}}_{j}}=	
	\sum_{h=1,\ h\not= z}^{H_{\hat\ell}} w^{{\hat\ell}+1}_{{\tilde c}h} g^\prime\biggl( a_{h}^{\hat\ell}(x^{p},\theta^{{\hat\ell}})\biggr)	 
	\frac{\partial  a_{h}^{{\hat\ell}}(x^{p},\theta^{{\hat\ell}})}{\partial   \sigma^{{\hat\ell}}_{j}}\\
	\nonumber 
	&&	\qquad\qquad\qquad\qquad\qquad+\lambda_0  w^{{\hat\ell}+1}_{{\tilde c}z}\; g^\prime\biggl(\hat a_{z}^{\hat\ell}(x^{p},\hat\theta^{{\hat\ell}})\biggr)	 
	\frac{\partial a_{z}^{{\hat\ell}}(x^{p},\theta^{{\hat\ell}})}{\partial   \sigma^{{\ell}}_{j}}
	\\ [0.3truecm]
	\nonumber
	\nonumber 
	&&	\qquad\qquad\qquad\qquad\qquad+\sum_{\hhat=H_{\hat\ell}+1}^{H_{\hat\ell}+K}\lambda_{\hhat-H_{\hat\ell}}  w^{{\hat\ell}+1}_{{\tilde c}z}\; g^\prime\biggl(\hat a_{z}^{\hat\ell}(x^{p},\hat\theta^{{\hat\ell}})\biggr)	 
	\frac{\partial  a_{z}^{{\hat\ell}}(x^{p},\theta^{{\hat\ell}})}{\partial   \sigma^{{\ell}}_{j}}
	\\ [0.3truecm]
	\nonumber
	&&	\qquad\qquad\qquad\qquad\qquad=\frac{\partial  a^{{\hat\ell}+1}_{\tilde c}(x^{p}, \theta^{{\hat\ell}+1}))}{\partial  \sigma^{{\hat\ell}}_{j}},
	\\ [0.3truecm]
	\nonumber
	&&\frac{\partial \hat a^{{\hat\ell}+1}_{\tilde c}(x^{p}, \hat\theta^{{\hat\ell}+1}))}{\partial \hat w^{\hat\ell}_{j i}}=	
	\sum_{h=1}^{H_{\hat\ell}} w^{{\hat\ell}+1}_{{\tilde c}h} g^\prime\biggl( a_{h}^{\hat\ell}(x^{p},\theta^{{\hat\ell}})\biggr)	 
	\frac{\partial  a_{h}^{{\hat\ell}}(x^{p},\theta^{{\hat\ell}})}{\partial  w^{\hat\ell}_{j i}}\\
	\nonumber
	&&	\qquad\qquad\qquad\qquad\qquad+\lambda_0  w^{{\hat\ell}+1}_{{\tilde c}z}\; g^\prime\biggl(\hat a_{z}^{\hat\ell}(x^{p},\hat\theta^{{\hat\ell}})\biggr)	 
\frac{\partial a_{z}^{{\hat\ell}}(x^{p},\theta^{{\hat\ell}})}{\partial  w^{{\ell}}_{jw}}
\\ [0.3truecm]
\nonumber 
&&	\qquad\qquad\qquad\qquad\qquad+\sum_{\hhat=H_{\hat\ell}+1}^{H_{\hat\ell}+{K}}\lambda_{\hhat-H_{\hat\ell}}  w^{{\hat\ell}+1}_{{\tilde c}z}\; g^\prime\biggl(\hat a_{z}^{\hat\ell}(x^{p},\hat\theta^{{\hat\ell}})\biggr)	 
\frac{\partial  a_{z}^{{\hat\ell}}(x^{p},\theta^{{\hat\ell}})}{\partial  w^{{\ell}}_{jw}}
\\ [0.3truecm]
\nonumber
	&&	\qquad\qquad\qquad\qquad\qquad=\frac{\partial  a^{{\hat\ell}+1}_{\tilde c}(x^{p}, \theta^{{\hat\ell}+1}))}{\partial w^{{\hat\ell}}_{ji}}.
	\nonumber
\end{eqnarray}

From the previous equalities, we get again

\begin{eqnarray}
	\label{eq46}
	&&\hspace*{-1.4cm}\frac{\partial \hat R_{emp}(\hat\theta^L)}{\partial \hat \sigma^{{\ell}}_{j}}  = \frac{1}{P} \sum_{p=1}^{P} \sum_{r=1}^{m} { \mathcal{L}^\prime(y^{p},  f_h(x^{p}, \theta^L))}
	\frac{\partial \hat a^{L}_r(x^{p}, \theta^L))}{\partial   \sigma^{{\ell}}_{j}}=\frac{\partial  R_{emp}(\theta^L)}{\partial   \sigma^{{\ell}}_{j}}  =0,\\ [0.3truecm]
	\label{eq47}
	&&\hspace*{-1.4cm}\frac{\partial \hat R_{emp}(\hat \theta^L)}{\partial \hat w^{{\ell}}_{ji}}  = \frac{1}{P} \sum_{p=1}^{P} \sum_{r=1}^{m}{ \mathcal{L}^\prime(y^{p},  f_r(x^{p}, \theta^L))} 
	\frac{\partial  a^{L}_r(x^{p},  \theta^L))}{\partial  w^{{\ell}}_{ji}}=\frac{\partial  R_{emp}( \theta^L)}{\partial  w^{{\ell}}_{ji}}  =0.
\end{eqnarray}
Then, the proof of the proposition follows from (\ref{eq34})-(\ref{eq47}).
$\hfill\square$

\section{Adding neurons to different layers}\label{sec:4}
Let us recall the definition of the mappings $\alpha$, $\beta$, and $\gamma$ given in section \ref{sec:3}, where we explicitly denote with $K_\ell$ the neurons added to the $\ell$-th layer, i.e.

\begin{eqnarray*}
 &&\hat\theta=\alpha_{\zeta {v}}(\theta;\ell,K_\ell), \\
	&&\hat\theta=\beta_{\zeta {s}}(\theta;\ell,K_\ell),\\
	&&\hat\theta=\gamma_{\lambda }(\theta;\ell,K_\ell).
\end{eqnarray*}

Let us consider $R= \{r_1,\dots,r_t\}\subseteq\{1,\dots,L\}$ and the finite set 
\[
\Gamma = \{K_{r_1},\dots,K_{r_t}\}.
\]
Therefore, we can define the composition of mappings
\begin{eqnarray*}
\theta_1 & = & \alpha_{\zeta {v}}(\theta\phantom{_1};r_1,K_{r_1}),\\
\theta_2 & = & \alpha_{\zeta {v}}(\theta_1;r_2,K_{r_2}),\\
\vdots   &   & \qquad\vdots\\
\theta_t & = & \alpha_{\zeta {v}}(\theta_{t-1};r_t,K_{r_t}),
\end{eqnarray*}
and call them $A_{\zeta v}(\theta,R,\Gamma)$, i.e.
\[
 \theta_t = \alpha_{\zeta {v}}\Big(\alpha_{\zeta {v}}\big(\dots \alpha_{\zeta {v}}(\theta;r_1,K_{r_1})\dots
 ;r_{t-1},K_{r_{t-1}}\big) ;r_t,K_{r_t}\Big) = A_{\zeta v}(\theta,R,\Gamma),
\]
which is a composition of mappings that, given $\theta$, produces $\hat\theta = \theta_t$. Analogously, for the $\beta$ and $\gamma$ mappings, we can define
\begin{eqnarray*}
\theta_1 & = & \beta_{\zeta {s}}(\theta\phantom{_1};r_1,K_{r_1}),\quad\ \ \theta_1   =   \gamma_{\lambda}(\theta\phantom{_1};r_1,K_{r_1}),\\
\theta_2 & = & \beta_{\zeta {s}}(\theta_1;r_2,K_{r_2}),\quad \ \
\theta_2   =   \gamma_{\lambda}(\theta_1;r_2,K_{r_2}),\\
\vdots   &   & \qquad\vdots\\
\theta_t & = & \beta_{\zeta {s}}(\theta_{t-1};r_t,K_{r_t}),\quad
\theta_t   =   \gamma_{\lambda}(\theta_{t-1};r_t,K_{r_t}),
\end{eqnarray*}
and call them $B_{\zeta s}(\theta,R,\Gamma)$ and $G_\lambda(\theta,R,\Gamma)$, i.e.
\[
 \theta_t = \beta_{\zeta {s}}\Big(\beta_{\zeta {s}}\big(\dots \beta_{\zeta {s}}(\theta;r_1,K_{r_1})\dots
 ;r_{t-1},K_{r_{t-1}}\big) ;r_t,K_{r_t}\Big) = B_{\zeta s}(\theta,R,\Gamma),
\]
and
\[
 \theta_t = \gamma_{\lambda}\Big(\gamma_{\lambda}\big(\dots \gamma_{\lambda}(\theta;r_1,K_{r_1})\dots
 ;r_{t-1},K_{r_{t-1}}\big) ;r_t,K_{r_t}\Big) = G_{\gamma}(\theta,R,\Gamma).
\]
\begin{proposition}\label{Prop1generale}
For every  point  $\theta$, it results 
\begin{eqnarray*}\label{equ-alfa}
&& \hat R_{emp}(\hat\theta)	=\hat R_{emp}(A_{\zeta v}(\theta,R,\Gamma))= R_{emp}(\theta),\\
\label{equ-beta}
&& \hat R_{emp}(\hat\theta)	=\hat R_{emp}(B_{\zeta {s}}(\theta,R,\Gamma) )= R_{emp}(\theta),\\
\label{equ-gamma}
&& \hat R_{emp}(\hat\theta)	=\hat R_{emp}(G_{\lambda }(\theta,R,\Gamma))= R_{emp}(\theta).
\end{eqnarray*}
\end{proposition}
\par\medskip\noindent
{\bf Proof}.
The proof follows by recursively applying the reasoning of Proposition \ref{Prop1}. $\hfill\Box$
\par\medskip\noindent

\begin{proposition}\label{Prop2generale}
	Let the point  $\theta$ be such that
\[
\nabla_\theta R_{emp}(\theta)=0,
\]
and let the point $\hat\theta$ be given by 
\[
	\hat\theta =B_{\zeta {0}}(\theta,R,\Gamma),
\]
or  
\[
	\hat\theta =G_{\lambda}(\theta,R,\Gamma).
\]
Then, it results: 
	\[
		 \nabla_{\hat\theta} \hat R_{emp}(\hat\theta)=0.
	\]
\end{proposition}
{\bf Proof}.
The proof follows quite easily by recursively applying Proposition \ref{Prop3}.$\hfill\Box$

\begin{remark}
Note that an analogous result to that of Proposition \ref{Prop2generale} can be obtained when different  mappings (i.e. $\beta$ or $\gamma$) are used when adding neurons to different layers of the network. In particular, we can define
\begin{eqnarray*}
\theta_1 & = & \xi^{(1)}(\theta\phantom{_1};r_1,K_{r_1}),\\
\theta_2 & = & \xi^{(2)}(\theta_1;r_2,K_{r_2}),\\
\vdots   &   & \qquad\vdots\\
\theta_t & = & \xi^{(t)}(\theta_{t-1};r_t,K_{r_t}),
\end{eqnarray*}
where each $\xi^{(i)}$ is either $\beta_{\zeta 0}$ or $\gamma_{\lambda}$, for $i=1,\dots,t$.
% \red{non è che vale anche se concateniamo $\alpha$ con $\beta$ e $\gamma$? Forse basterebbe dire che si può fare?}
\end{remark}

Considering again the results of Proposition \ref{Prop2generale}, it can be noticed that the number of the manifolds of useless stationary points in a given network grows exponentially with the network dimension.

\section{The incremental training algorithm}\label{sec:5}

In this section, we formally state our proposed incremental training algorithm \ref{alg:CMA2}. The results described in the previous section have proved that every stationary point of a smaller network corresponds to certain stationary points of a larger network. In particular we have shown that, given a network with $L$ layers, it is possible to add a given number of neurons on some of (or all) the layers preserving stationarity. Indeed, if $\theta_*$ is stationary, then $\hat\theta_* = B_{\zeta 0}(\theta_*,R,\Gamma)$ or $\hat\theta_* = G_{\lambda}(\theta_*,R,\Gamma)$ are stationary in the bigger network. However, the latter may not correspond to a global minimum. 

For the sake of simplicity, in this section we consider neural networks with a single hidden layer. In this case, we have $L=1$ hidden layer with $H$ neurons and let
\[
  R = \{1\},\quad\mbox{and}\quad \Gamma = \{K\},
\]
so that the mappings are
\begin{eqnarray*}
\theta^{(H+K)} & = & \alpha_{\zeta v}(\theta^{(H)};1,K) = \alpha_{\zeta v}(\theta^{(H)};K),\\
\theta^{(H+K)} & = & \beta_{\zeta s}(\theta^{(H)};1,K) = \beta_{\zeta s}(\theta^{(H)};K),\\
\theta^{(H+K)} & = & \gamma_{\lambda}(\theta^{(H)};1,K) = \gamma_{\lambda}(\theta^{(H)};K),
\end{eqnarray*}
where we denoted with $\theta^{(H)}$ the vector of parameters of the network with $H$ neurons.
The idea behind our Incremental Training Algorithm (ITA) is training a network of a given dimension starting from a smaller network and progressively increasing the number of neurons in the hidden layer. At every iteration, the training of the larger network is performed by properly choosing both the starting point and the minimization technique. Regarding the starting points, we use the mapping $\alpha_{\zeta v}$, which guarantees that the objective function value doesn't change. Instead, the gradient is unlikely to be equal to zero. Therefore, any descent algorithm that is able to produce a sequence $\{\theta_h^{(H+K)}\}$ such that the objective function value satisfies
$$R_{emp}(\theta_{h}^{(H+K)}) < R_{emp}(\theta_{0}^{(H+K)}) =R_{emp}(\alpha_{\zeta v}(\theta^{(H)};K)).$$ 
is not attracted by those useless stationary points described in Proposition \ref{Prop2generale}.
Considering that 
\[
R_{emp}(\alpha_{\zeta v}(\theta^{(H)}_*;K)) = 
R_{emp}(\beta_{\zeta 0}(\theta^{(H)}_*;K)) = 
R_{emp}(\gamma_{\lambda}(\theta^{(H)}_*;K))
%
%R_{emp}(\alpha_{\tilde w}(\theta_{*}^{(H)})) = 
%R_{emp}(\beta_{(0,\xi)}(\theta_{*}^{(H)}))) = 
%R_{emp}(\gamma_{\lambda}(\theta_{*}^{(H)}))),
\]
the algorithm has the ability to escape from the stationary points generated  by the mappings $\beta_{\zeta 0}(\theta^{(H)}_*;K)$ and $\gamma_{\lambda}(\theta^{(H)}_*;K)$.  

\begin{algorithm}[htb]
\begin{algorithmic}[1]
	\State {\bf Data}: $H_0, H_{\max} \in\mathbb{N}$, $\{K_k\}\subset\mathbb{N}$, $\{\tau_k\}$.
	\State Set $\theta^{(H_0)}\propto U(0,1)$, $k\leftarrow 0$.
	\While {$H_k \leq H_{\max}$}
	\State Compute 
	\[
	\bar\theta^{(H_k)}\ \mbox{such that}\ \left\|\nabla R_{emp}\left(\bar\theta^{(H_k)}\right)\right\| \leq \tau_k.
	\]	
	\State Let $H_{k+1} \leftarrow \min\{H_k + K_k,H_{\max}\}$.
%	\State Choose $\tilde w_i \propto U(0,1)$, for  $i=1,\dots,\min\{H_{\max}-H_k,K_k\}$
 	\State Choose $(\zeta_i,v_{i}) \propto U(0,1)$, for  $i=1,\dots,\min\{H_{\max}-H_k,K_k\}$
	\State set $\theta^{(H_{k+1})} \leftarrow  \alpha_{\zeta v}(\bar\theta^{(H_k)})$, s.t.
	\[
	   R_{emp}(\theta^{(H_{k+1})}) = R_{emp}(\bar\theta^{(H_k)}),\ \left\|\nabla R_{emp}\left(\theta^{(H_{k+1})}\right)\right\| > \tau_k
	\]
	\State Set $k\leftarrow k+1$.
	\EndWhile
\end{algorithmic}
\caption{Incremental Training Algorithm (ITA)}
\label{alg:CMA2}
\end{algorithm}

\section{Numerical experiments}\label{sec:6}
In this section, we report numerical results to support the observations made above on our incremental training algorithm (ITA). To verify its scalability, we test the method on a set of standard test problems of different sizes mostly taken from the OpenML and UCI Machine Learning repositories as reported in Table \ref{table:5.1}, with the exception of the Power Consumption dataset, which represents the hourly power consumption of one of the largest energy provider in Italy.
\begin{table}[ht]\small
\caption{List of the test problems considered} 
\centering 
\begin{tabular}{l l c c} 
\hline\hline 
Name & Type & \# Instances & \# Attributes\\ [0.5ex]
\hline
Adult & Classification & 48842  & 14 \\
Ailerons & Regression & 13750 & 40 \\
Appliances Energy Prediction & Regression & 19735 & 29 \\
Arcene & Classification & 200 & 10000 \\
BlogFeedback & Regression & 60021 & 281 \\
Boston House Prices & Regression & 21613 & 19 \\
Breast Cancer Wisconsin (Diagnostic) & Classification & 569 & 32 \\
CIFAR 10 & Classification & 20000 & 3072 \\
Gisette & Classification & 13500 & 5000 \\
Iris & Classification & 150 & 4 \\
MNIST Handwritten Digit & Classification & 70000 & 784 \\
Mv & Regression & 40768 & 10 \\
QSAR Oral Toxicity & Classification & 8992 & 1024 \\
Power Consumption & Regression & 4520 & 347 \\
YearPred & Regression & 515345 & 90 \\ [1ex] 
\hline 
\end{tabular}
\label{table:5.1}
\end{table}
\par\smallskip
Our incremental training algorithm was implemented in Python (v.3.8.5) using a XPS 15 7590 Intel(R) Core(TM) i7-9750H CPU, 2.60GHz, 16 GB RAM. 
For each test problem, we assessed both the standard and incremental approaches by training a neural network 
using the Pytorch (v.1.7.0) library with $\tanh(t)$ as activation function for each neuron in the hidden layer, a linear activation function for the output layer and the mean squared error as loss function. For both the methods considered, we used the L-BFGS optimization algorithm available in Pytorch with a maximum number of training epochs {\tt maxit} $=1000$. In particular, for the standard case we adopted the following parameters configuration:
\begin{itemize}
    \item[-] $H = 100$ neurons in the hidden layer;
    \item[-] a tolerance {\tt tol} $=10^{-6}$ in the stopping condition, i.e.
    \[
       \|\nabla R_{emp}(\theta^{(H)})\|_{\infty} \leq 10^{-6}.
    \]
\end{itemize}

 %and on the decrease of the gradient norm equal to $\|\nabla f(x_{k})\|_{\infty} \leq 10^{-6}$.  
 
Instead, for the incremental method (ITA) we made the following choices:
 \begin{itemize}
 \item[-] initial number of neurons in the hidden layer $H_{0} = 10$ and $K_k = H_k$;
 \item[-] stopping criterion for the intermediate networks
  \[
   \|\nabla R_{emp}(\theta^{(H_k)}_{h})\|_{\infty} \leq 10^{-1}\|\nabla R_{emp}(\theta^{(H_k)}_{0})\|_{\infty}.
 \]
\end{itemize}
From numerical experience, we noticed that a better performance can be obtained by also adding the following criterion for the intermediate networks: 
$${|R_{emp}(\theta^{(H_k)}_{h}) - R_{emp}(\theta^{(H_k)}_{h-1})|} \leq 10^{-2}.$$
In the following paragraphs, we detailed the numerical results obtained from both the approaches considered, by firstly reporting their whole performance profiles in order to give a global perspective of the methods. Then, we gave a more detailed description of the proposed methods by reporting the training loss and box plot distribution for each single test problem evaluated.

\paragraph{Performance profiles}
To compare the two variants of training methods, we adopted the performance profiles proposed in \cite{dolan:2002}. 
Let $\cal P$ be a set of $n_p$ problems and $\cal S$ be a set of $n_s$ solvers that can be used to solve problems in $\cal P$. For each $s\in\cal S$ and $p\in\cal P$, let $t_{ps}$ denote the performance index (i.e. the final function value obtained by the  solver). The performance ratio is then defined as
\[
  r_{ps} = \frac{t_{ps}}{\min_{s\in\cal S}\{t_{ps}\}}.
\]
The performance profile relative to solver $s$ is defined as
\[
 \rho_s(\alpha) = \frac{1}{|\cal P|}\Big|\Big\{p\in{\cal P}:\ r_{ps} \leq \alpha\Big\}\Big|,
\]
where $\alpha\ge 1$. Basically, for each solver a performance profile reports the
percentage of problems for which a final function value is obtained which is within $\alpha$-times the function value attained by the best solver. Hence, the uppermost curve in the profiles denotes better performances of the corresponding algorithm. In our experimentation, the solvers are the incremental and standard methods whereas the problem set is composed of 10 replica for each of the 15 test problems, thus amounting to 150 problems. We reported the performance profiles of the two methods in Fig. \ref{fig2}, \ref{fig3}, \ref{fig4} respectively with {\tt maxit} = $200, 500, 1000$. All these plots seem to suggest that our incremental training algorithm is significantly more efficient than the standard approach from a global behaviour perspective on all the datasets assessed. 
\paragraph{Boxplot and training loss}
To further understand the difference between the two approaches, we plotted the training loss and the boxplot distribution over 10 replica of each test problem assessed. Fig. \ref{fig5}  describes the performance of the two approaches in terms of the boxplot distribution with {\tt maxit} = $1000$. In this case, we can notice that in some datasets, i.e. Arcene, Iris, Qsar and Power Consumption, our incremental method is definitely better that the standard approach. Instead, for the Mv and YearPred cases, our approach performs worse than the standard one. In all the other cases, the performances of the two approaches are comparable. However, it is important to take into account that the information contained in the boxplots denotes an ensemble behaviour hiding the dynamics of the performance of the epochs, as we can appreciate from the loss function progress in
Fig. \ref{fig6}, which describes  the performance of the two approaches in terms of the training loss with {\tt maxit} = $1000$.
After some initial iterations, following the incremental method consistently offers an additional progress, allowing to reduce the training loss more rapidly than the standard case. Even in those cases where it performs worse than the standard method, the incremental approach allows to reach the same training value of the standard method in fewer iterations, as shown in Fig. \ref{fig7}, \ref{fig8}. These preliminary results seem to suggest that the incremental approach may be useful to avoid the early stagnation experienced by the standard method.

\section{Conclusions}\label{sec:7}
In this paper we extended the result of \cite{fukumizu:2000} by characterizing the structure of undesirable stationary points to deep multilayer neural networks, i.e. networks with more than one layer. More precisely, we show that the structure of such manifolds of stationary points implies that their number grows exponentially with the dimension of the network.

Moreover, a novel incremental approach that avoids such undesirable stationary points is proposed. Unlike the traditional method, the main advantage of the proposed scheme is training a network of a given dimension starting from a smaller network and progressively increasing the number of neurons. Therefore, this method is able to escape from useless stationary points of a certain neural network by progressively exploit the information contained in the smaller networks. Numerical experiments on a significant number of test problems show the good performances of the proposed method when compared with a standard approach. In particular, our incremental scheme seems to be able to avoid the early stagnation experienced by the standard method on all the datasets assessed, thus suggesting more efficiency from a global behaviour perspective. However, the numerical results obtained refer to the case of a network with only one single hidden layer. The impact of the proposed method in deep multilayer neural networks training has not been studied yet in terms of numerical experience and will be developed in a future study.

\bibliographystyle{plain}
\bibliography{incremental.bib}

\appendix
\section*{Appendix}

\begin{figure}[H]
    \begin{center}
    \caption{Performance profiles of the standard and incremental training algorithms over 100 maximum iterations}
    \includegraphics[width=\textwidth]{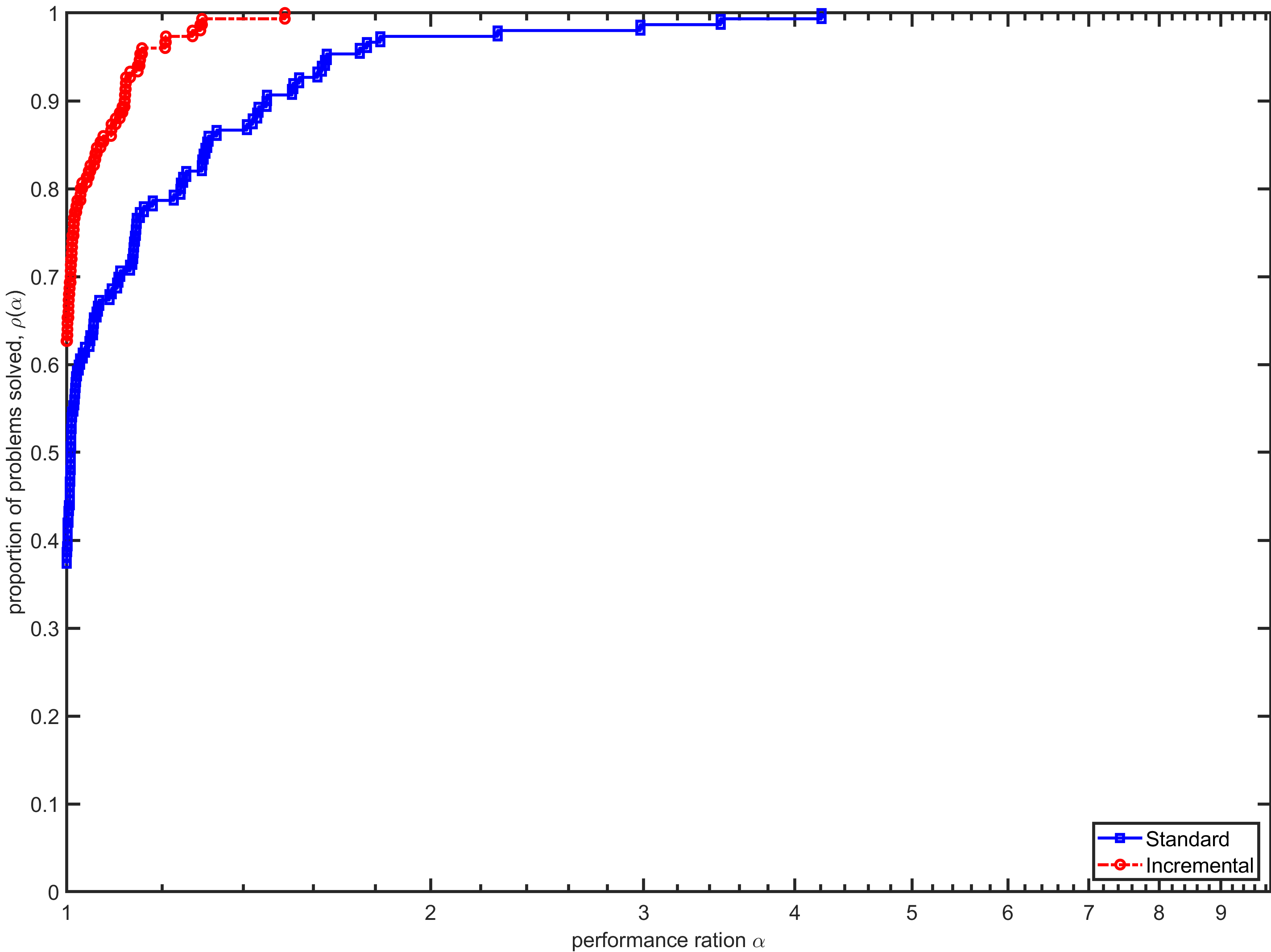}
    \label{fig2}
    \end{center}
\end{figure}

\begin{figure}[H]
    \begin{center}
    \caption{Performance profiles of the standard and incremental training algorithms over 500 maximum iterations}
    \includegraphics[width=\textwidth]{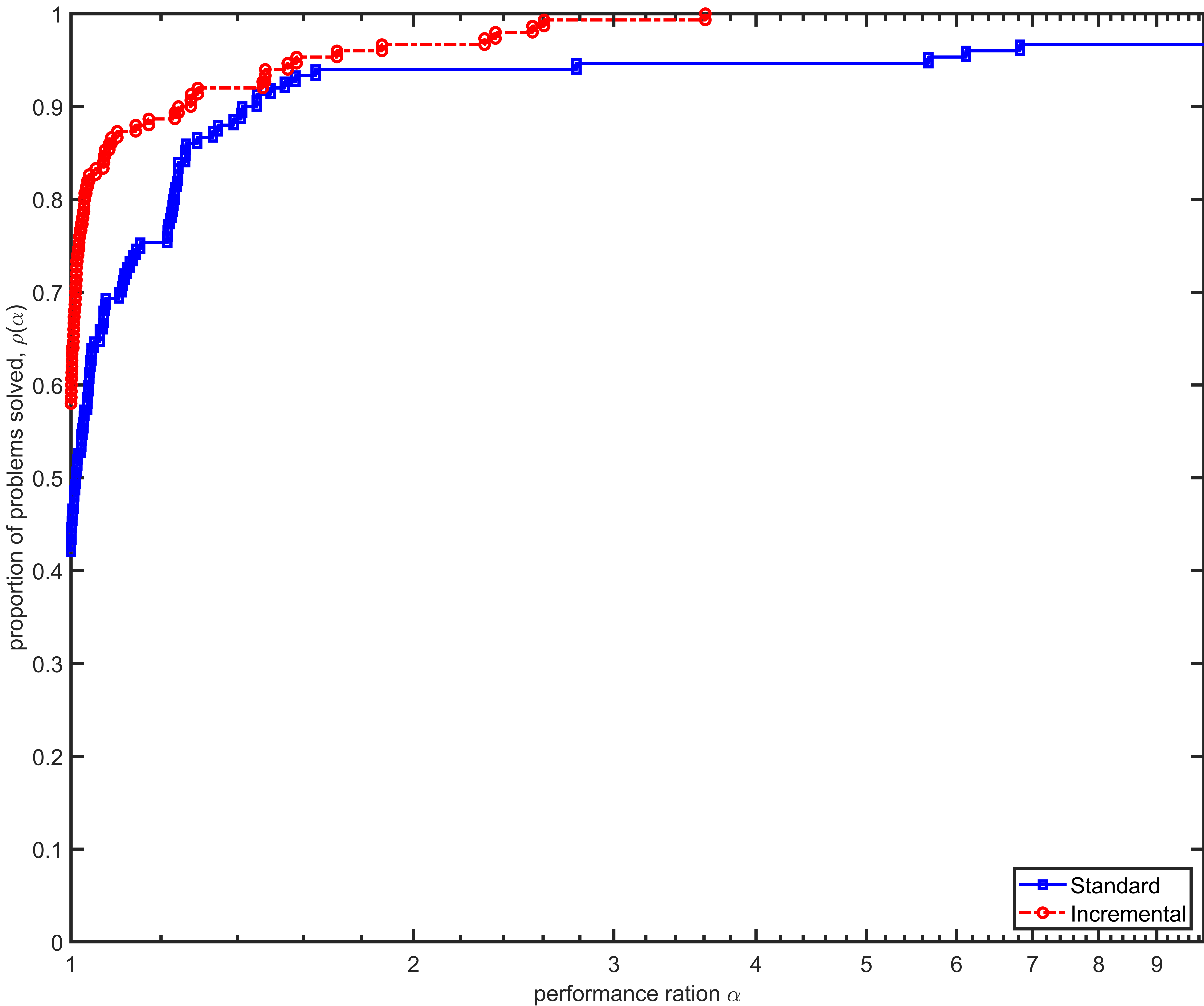}
    \label{fig3}
    \end{center}
\end{figure}

\begin{figure}[H]
    \begin{center}
    \caption{Performance profiles of the standard and incremental training algorithms over 1000 maximum iterations}
    \includegraphics[width=\textwidth]{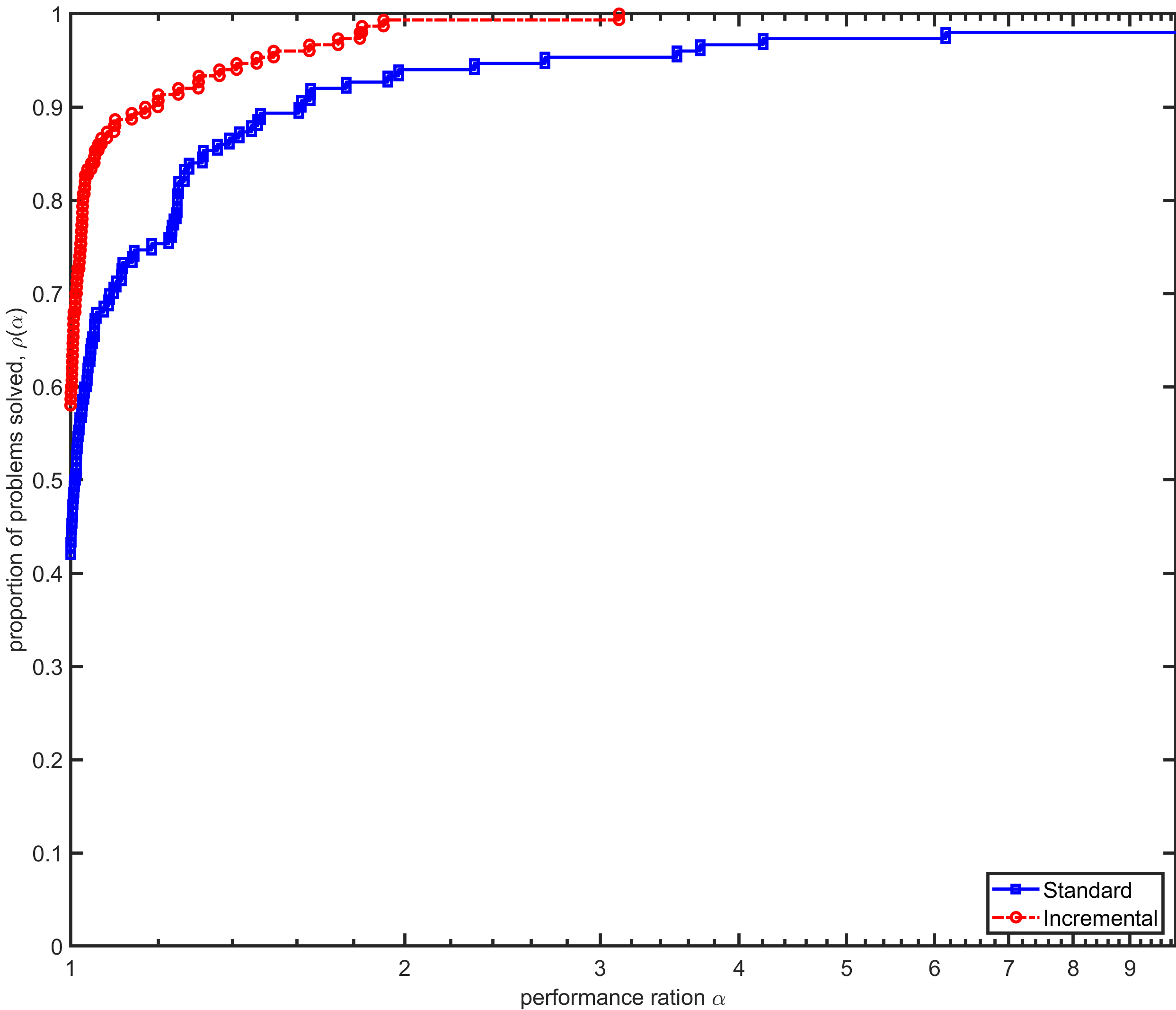}
    \label{fig4}
    \end{center}
\end{figure}

\begin{figure}[H]
    \centering
    \caption{Boxplot distribution for the standard and incremental training algorithms over 1000 maximum iterations}
    \advance\leftskip-2.5cm
    \includegraphics[width=1.4\textwidth]{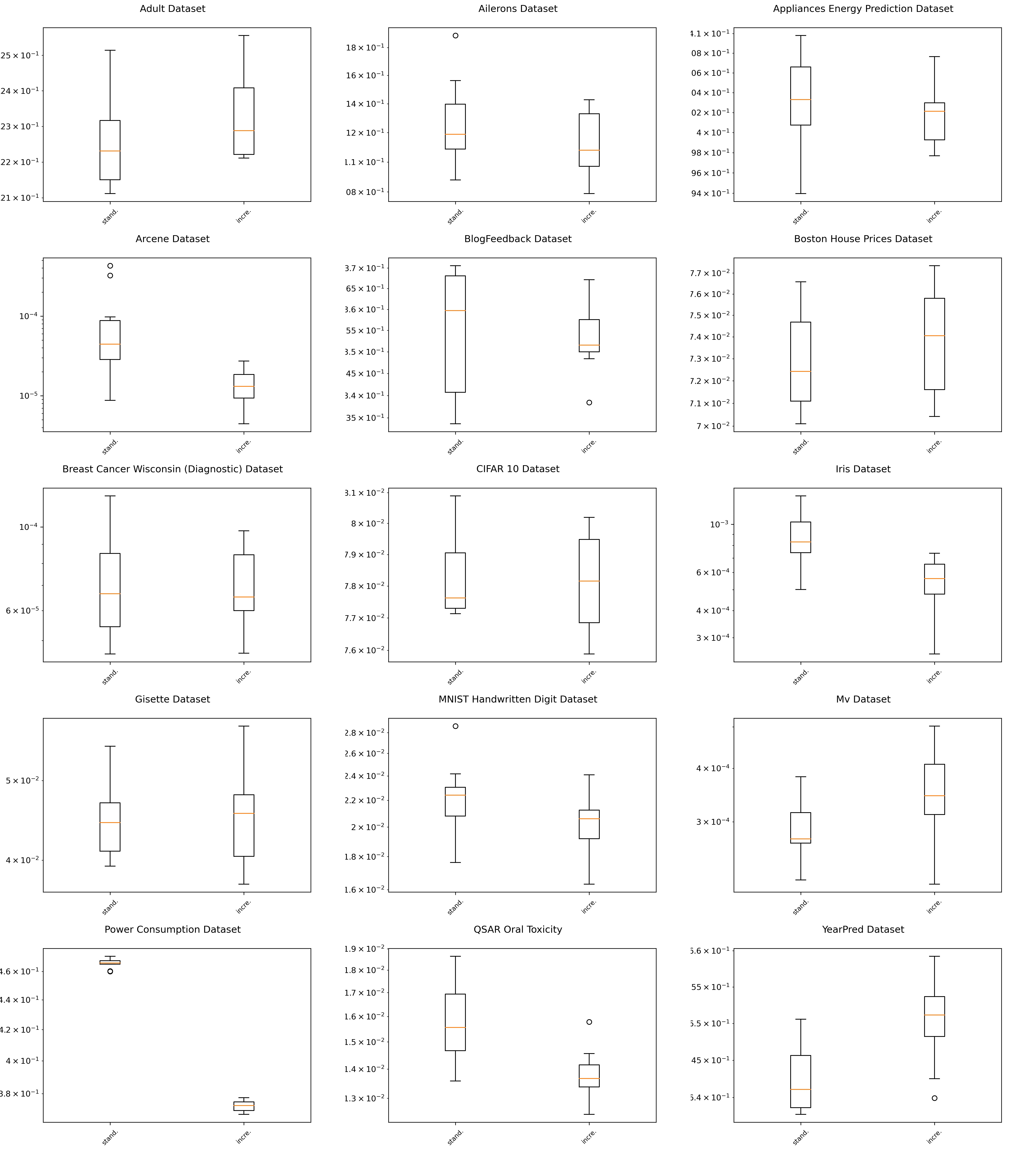}
    \label{fig5}
\end{figure}

\begin{figure}[H]
    \centering
    \caption{Training loss for the standard and incremental training algorithms over 1000 maximum iterations}
    \advance\leftskip-2.5cm
    \includegraphics[width=1.4\textwidth]{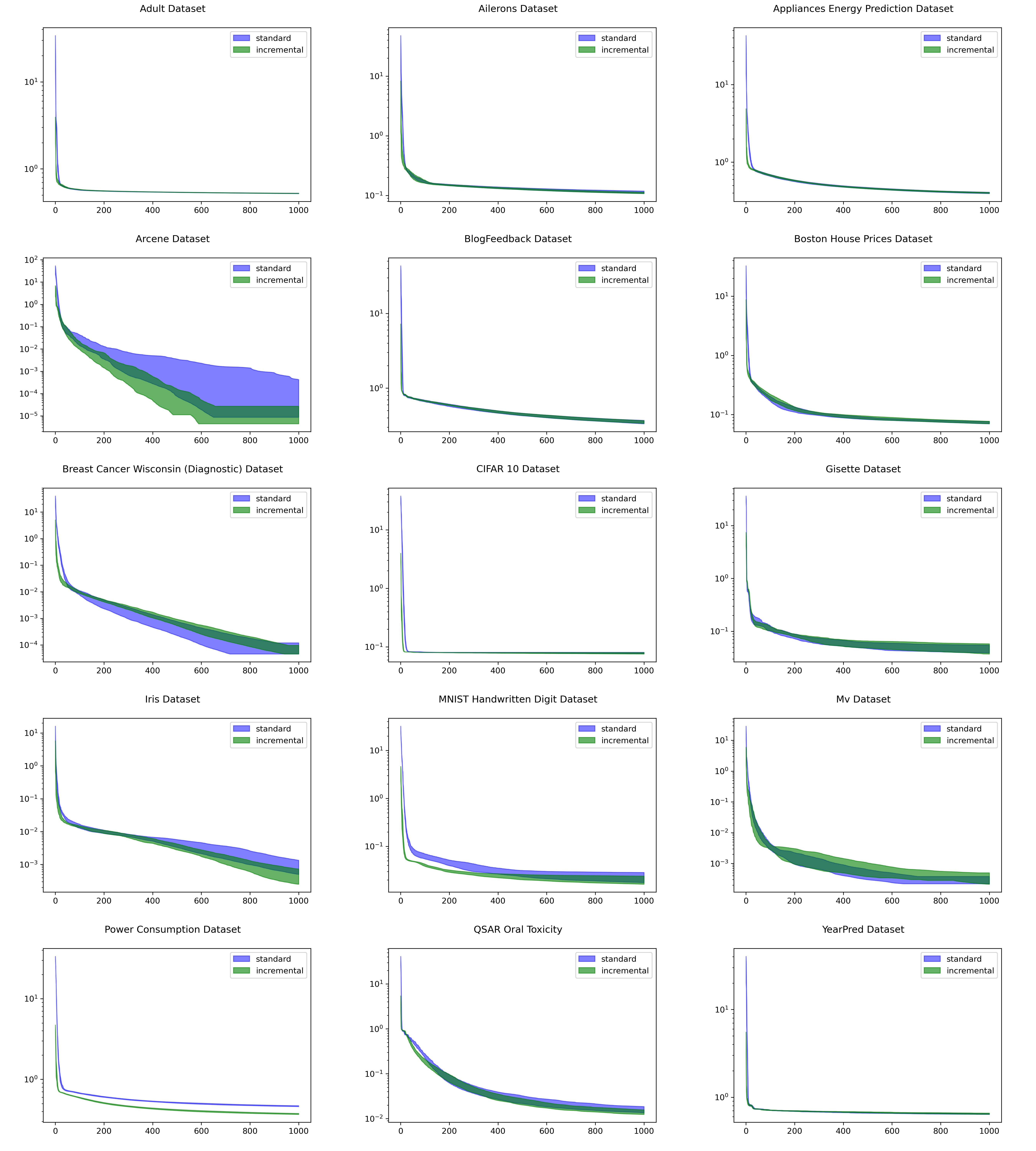}
    \label{fig6}
\end{figure}

\begin{figure}[H]
    \centering
    \caption{Training loss for the standard and incremental training algorithms over 100 maximum iterations}
    \advance\leftskip-2.5cm
    \includegraphics[width=1.4\textwidth]{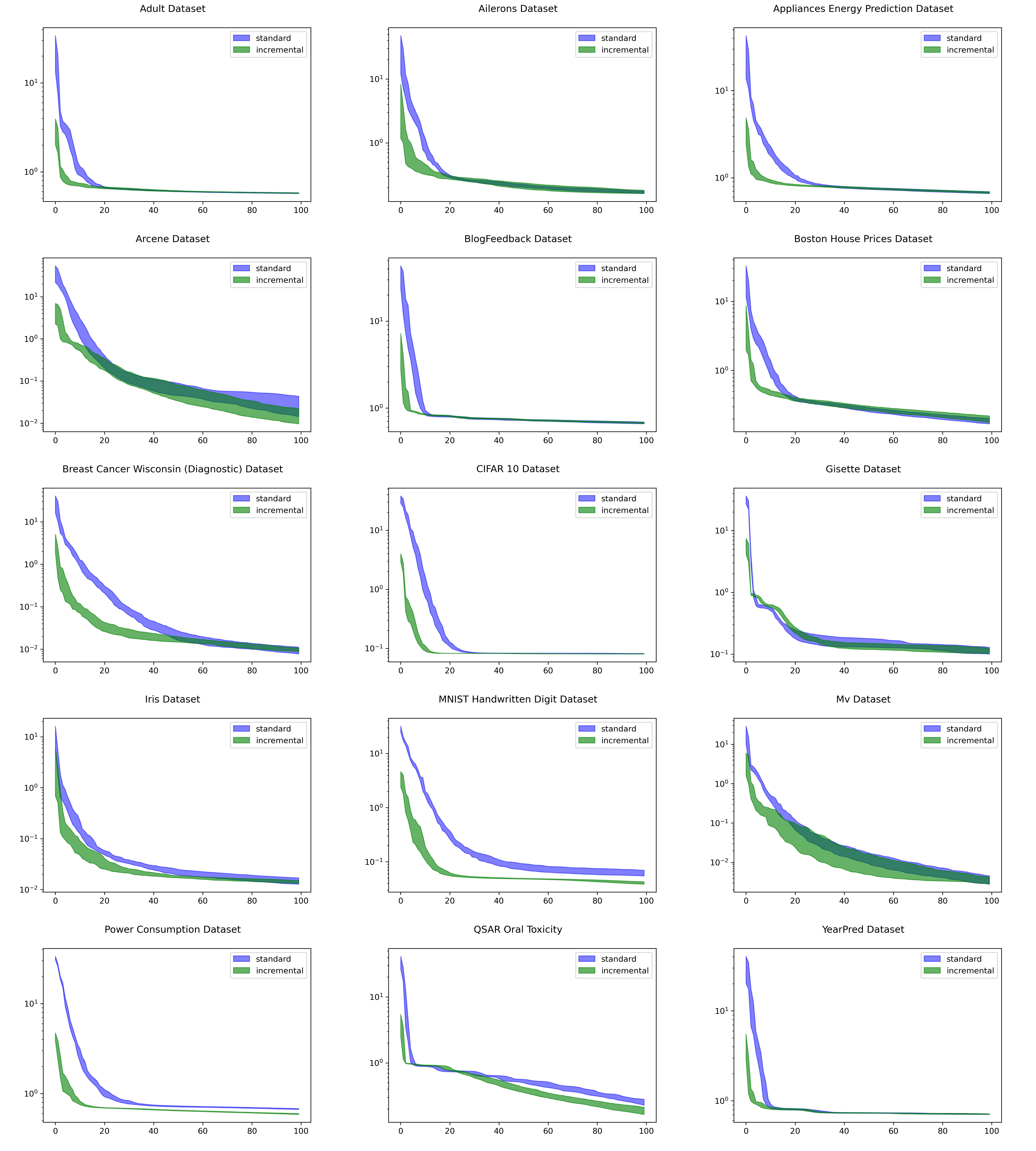}
    \label{fig7}
\end{figure}

\begin{figure}[H]
    \centering
    \caption{Training loss for the standard and incremental training algorithms over 500 maximum iterations}
    \advance\leftskip-2.5cm
    \includegraphics[width=1.4\textwidth]{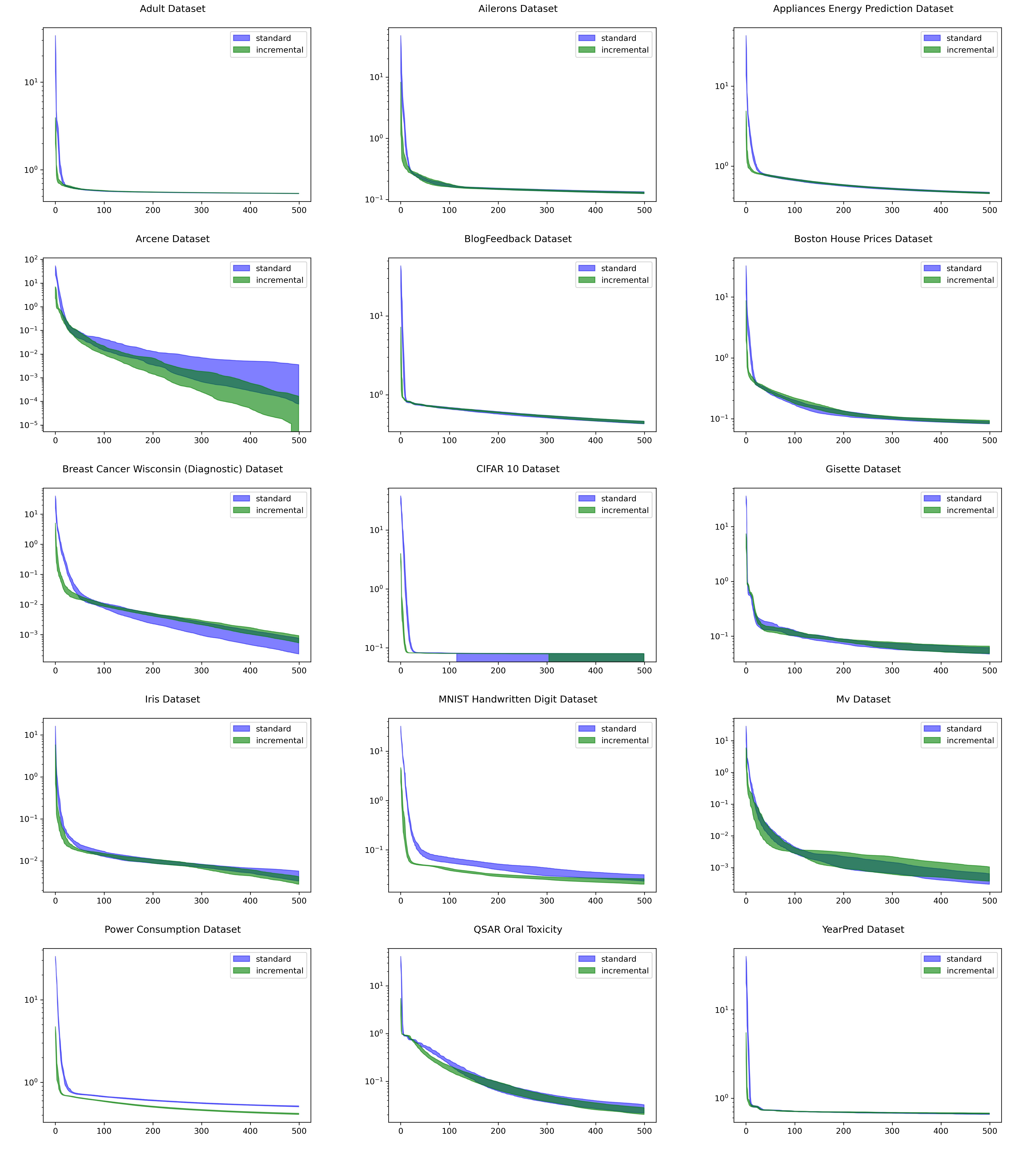}
    \label{fig8}
\end{figure}

\end{document}